\documentclass{article}
\usepackage{amssymb}
\usepackage{amsfonts}
\usepackage{amsmath}
\usepackage[a4paper,ignoreheadfoot,noheadfoot,nomarginpar]{geometry}

\newtheorem{theorem}{Theorem}

\newtheorem{corollary}[theorem]{Corollary}

\newtheorem{example}[theorem]{Example}

\newtheorem{lemma}[theorem]{Lemma}

\newtheorem{proposition}[theorem]{Proposition}

\newenvironment{proof}[1][Proof]{\noindent\textbf{#1.} }{\ \rule{0.5em}{0.5em}}
\input{tcilatex}
\begin{document}

\title{Polynomial Weingarten Tubular Surfaces in Euclidean, Hyperbolic and
Lorentzian 3-spaces}
\author{A.P.Barreto and F.Gasparotto \\
Universidade Federal de S\~{a}o Carlos}
\maketitle

\begin{abstract}
In this article we fully classify regular tubular surfaces in Euclidean,
Lorentzian and hyperbolic 3-spaces whose Gaussian and mean curvatures $K$
and $H$ verify a polynomial relation. More precisely, we determine the set $%
\mathcal{S}\left( Q\right) $ of all regular tubular surfaces whose
curvatures verify a given polynomial relation $Q\left( K,H\right) =0$, and
the set $\mathcal{Q}\left( S\right) $ of all polynomial relations vanished
by the curvatures of a given tubular surface $S$.
\end{abstract}

\section{Introduction\label{Introduction}}

\qquad An important research problem in classical differential geometry is
to find relations between global geometric/topological properties of
immersed surfaces and hypothesis on its curvatures. In this direction, a
topic that has received much attention over the years is that of Weingarten
surfaces, that is, two dimensional surfaces in three dimensional spaces
whose Gaussian curvature $K$ and mean curvature $H$ (or its principal
curvatures $k_{1}\leq k_{2}$) satisfy a smooth non trivial relation%
\begin{equation*}
\Phi \left( K,H\right) \equiv 0\mathbb{\qquad }\text{(or }\Phi \left(
k_{1},k_{2}\right) \equiv 0\text{).}
\end{equation*}%
When $\Phi $ is a polynomial in $\mathbb{R}\left[ x,y\right] $, such
surfaces are called polynomial Weingarten surfaces. It is important to
observe that the class of polynomial Weingarten surfaces includes some
relevant and well studied families like minimal surfaces ($H\equiv 0$),
constant mean curvature surfaces ($H\equiv c\neq 0$), constant Gaussian
curvature surfaces ($K\equiv c$), linear Weingarten surfaces ($aK+bH\equiv c$%
) and surfaces with second fundamental form of constant length ($%
4H^{2}-2K\equiv c>0$).

\bigskip

Despite being an ancient topic, there is still much to be discovered about
Weingarten surfaces. Given a three dimensional space, the researchers in
this area are interested in (some of) the following problems:

\begin{enumerate}
\item[i.] \textbf{(}$\mathcal{S}\left( Q\right) $\textbf{\ problem)} Given a
smooth relation $\Phi \left( x,y\right) $, describe the set $\mathcal{S}%
\left( Q\right) $ of all regular surfaces whose curvatures verify $\Phi
\left( K,H\right) \equiv 0$.

\item[ii.] \textbf{(}$\mathcal{Q}\left( S\right) $\textbf{\ problem)} Given
a regular surface $S$, describe the set $\mathcal{Q}\left( S\right) $ of all
smooth relations $\Phi \left( x,y\right) $ verifying $\Phi \left( K,H\right)
\equiv 0$.

\item[iii.] \textbf{(Linear problem) }Determine if a given regular surface $%
S $ is a linear Weingarten surface, that is, if there is a non trivial
linear relation $\Phi \left( x,y\right) =ax+by+c$ in $\mathcal{Q}\left(
S\right) $.

\item[iv.] \textbf{(True non linear problem) }Determine if the curvatures of
a given regular surface $S$ vanish a true non trivial linear relation, that
is, if there is in $\mathcal{Q}\left( S\right) $ a smooth relation that
cannot be divided by another non constant relation in $\mathcal{Q}\left(
S\right) $.

\item[v.] \textbf{(Second fundamental form problem)} Determine the set of
all regular surfaces with second fundamental form of constant length. In
particular, determine if there are, in the Euclidean $3$-space, regular
surfaces of this type other than spheres and right cylinders.
\end{enumerate}

\bigskip

Due to its complexity, works about the problems (i) and (ii) are very rare.
In general, the results in the literature classify a particular familiy of
surfaces, in a classical environment and verifying one specific (linear in
general) relation. For instance, in \cite{BFH} the first author and its
colaborators classified rotational surfaces, in the Euclidean space, with
second fundamental form of constant length (for other works of this kind
see, for example, \cite{HW}, \cite{KS}, \cite{KT}, \cite{L2}, \cite{L3}, 
\cite{RE}, \cite{V}). In \cite{K}, \cite{DK} and \cite{DGW1}, however, the
authors were able to deal with the $\mathcal{S}\left( Q\right) $ problem for
ruled and translation Weingarten surfaces in Euclidean and Lorentzian $3$%
-spaces. Unfortunately, their techniques were not effective to treat the $%
\mathcal{Q}\left( S\right) $ problem.

\bigskip

In this paper we used an alternative algebraic approach that allowed us to
deal with the five problems described above for regular tubular surfaces
(see section \ref{Tubular Surface} for the definition). More precisely, we
fully classify regular tubular surfaces in Euclidean, Hyperbolic and
Lorentzian $3$-spaces whose Gaussian and mean curvatures verify an arbitrary
polynomial relation (see \cite{ASK}, \cite{KT}, \cite{KzT}, \cite{KYT}, \cite%
{L1}, \cite{S}, for a sample of works about tubular surfaces). In the
Euclidean $3$-space, our results read as follows:

\begin{theorem}[Euclidean $\mathcal{S}\left( Q\right) $ problem]
\label{TeoTubIntro}Given a polynomial $Q\left( x,y\right) \in \mathbb{R}%
\left[ x,y\right] $, denote by $\mathcal{S}_{E}\left( Q\right) $ the set of
all regular tubular surfaces in the Euclidean $3$-space verifying $Q\left(
K,H\right) \equiv 0$. Then, the elements of $\mathcal{S}_{E}\left( Q\right) $
are:

\begin{enumerate}
\item[i.] The right cylinders with radius in the set (called \textbf{%
Euclidean radius} of $Q$)%
\begin{equation*}
Rad_{E}\left( Q\right) =\left\{ r\in \left( 0,+\infty \right) \text{ };\text{
}Q\left( 0,\tfrac{1}{2r}\right) =0\right\} \text{;}
\end{equation*}

\item[ii.] The regular tubular surfaces with radius in the set (called 
\textbf{Euclidean star radius} of $Q$)%
\begin{equation*}
Rad_{E}^{\ast }\left( Q\right) =\left\{ r\in Rad_{E}\left( Q\right)
\,;\,Q\left( x,y\right) \in \left\langle xr^{2}-2ry+1\right\rangle \right\} 
\text{.}
\end{equation*}
\end{enumerate}
\end{theorem}

\begin{corollary}[Euclidean]
Consider a polynomial $Q\left( x,y\right) \in \mathbb{R}\left[ x,y\right] $.
There is a regular tubular surfaces in the Euclidean $3$-space verifying $%
Q\left( K,H\right) \equiv 0$ if and only if $Rad_{E}\left( Q\right) \neq
\emptyset $. Moreover, the radii of all such surfaces belongs to $%
Rad_{E}\left( Q\right) $.
\end{corollary}

\bigskip

Unlike the technique proposed by Kh\"{u}nel and Dillen in the articles
mentioned above, our Algebraic approach allows us to determine the set of
all polynomials in$\mathbb{\ \mathbb{R}}\left[ x,y\right] $ that vanishes
the Gaussian and the mean curvatures of a given regular tubular surface.

\begin{theorem}[Euclidean $\mathcal{Q}\left( S\right) $ problem]
\label{QSIntro}Given a regular tubular surface $S$ of radius $r>0$ in the
Euclidean $3$-space, denote by $\mathcal{Q}_{E}\left( S\right) $ the set of
all polynomials $Q\in \mathbb{R}\left[ x,y\right] $ verifying $Q\left(
K,H\right) \equiv 0$.

\begin{enumerate}
\item[i.] If $S$ is a right cylinder, then $\mathcal{Q}_{E}\left( S\right)
=\left\{ Q\in \mathbb{R}\left[ x,y\right] \,;\,r\in Rad_{E}\left( Q\right)
\right\} $;

\item[ii.] If $S$ is not a right cylinder, then $\mathcal{Q}_{E}\left(
S\right) $ is the ideal in $\mathbb{R}\left[ x,y\right] $ generated by the
polynomial $xr^{2}-2ry+1$.
\end{enumerate}
\end{theorem}

\begin{corollary}[Euclidean Linear problem]
Every regular tubular surface in Euclidean $3$-space is a linear Weingarten
surface. More precisely, their curvatures verify%
\begin{equation*}
Kr^{2}-2rH+1=0\text{,}
\end{equation*}%
where $r>0$ are their radii.
\end{corollary}

\begin{corollary}[Euclidean Nonlinear problem]
\label{Nonlinear intro}The right cylinders are the only tubular surfaces in
Euclidean $3$-space whose Gaussian and mean curvatures vanishes true
nonlinear polynomial relations. More precisely, if $S$ is a regular tubular
surface that is not a right cylinder and $Q$ is a polynomial in $\mathcal{Q}%
\left( S\right) $, then $Q$ has a non trivial divisor in $\mathcal{Q}\left(
S\right) $.
\end{corollary}

\bigskip

Many particular results can be obtained as consequence of our main theorems.
We feature here the classification of regular tubular surfaces verifying
linear relations and the classification of tubular surfaces with second
fundamental form of constant length. We address the reader to the section %
\ref{main result and applications} for more details and other results.

\begin{corollary}[Euclidean]
Given real numbers $a,b,c$ such that $\left( a,b\right) \neq \left(
0,0\right) $, consider the polynomial%
\begin{equation*}
Q\left( x,y\right) =ax+by-c\in \mathbb{R}\left[ x,y\right] \text{.}
\end{equation*}%
Then, $\mathcal{S}_{E}\left( Q\right) \neq \emptyset $ if and only if $b,c=0$
or $bc>0$. More precisely, for $\Delta =b^{2}+4ac$, we have:

\begin{enumerate}
\item[i.] If $b,c=0$, then $\mathcal{S}_{E}\left( Q\right) $ is the set of
all right cylinders of any radius;

\item[ii.] If $bc>0$ and $\Delta =0$, then $\mathcal{S}_{E}\left( Q\right) $
is the set of all tubular surfaces of radius $\frac{b}{2c}$;

\item[iii.] If $bc>0$ and $\Delta \neq 0$, then $\mathcal{S}_{E}\left(
Q\right) $ is the set of all right cylinders of radius $\frac{b}{2c}$.
\end{enumerate}
\end{corollary}

\begin{corollary}[Euclidean]
The right cylinders are the unique non degenerated regular tubular surfaces
in th Euclidean $3$-space with second fundamental form of constant length.
More precisely, the unique non degenerated regular tubular surfaces in the
Lorentzian $3$-space verifying $\left\vert \mathcal{A}\right\vert =c>0$ are
the right cylinders of radius $\frac{1}{c}$
\end{corollary}

\bigskip

The article is organized as follows: in section \ref{Curves and surfaces},
we review basic concepts about the Lorentzian and Hyperbolic spaces. In
section \ref{Tubular Surface} we present the definition of tubular surfaces,
we discuss the regularity of its parametrizations and we finish with the
calculation of its curvatures. Section \ref{Polynomials results} is a
technical section devoted to all polynomial results used in this paper.
Finally, in section \ref{main result and applications} we present our main
theorems and some of its applications.

\section{Curves and Surfaces in the Lorentzian and the Hyperbolic Spaces 
\label{Curves and surfaces}}

\qquad The \textbf{Lorentzian }$n$\textbf{-space (}$n\geq 2$), denoted here
by $\mathbb{L}^{n}$, is the semi-Riemannian manifold obtained endowing $%
\mathbb{R}^{n}$ with the bilinear form of index $1$ given by%
\begin{equation*}
\left\langle u,v\right\rangle _{1}=\sum_{i=1}^{n-1}u_{i}v_{i}-u_{n}v_{n},
\end{equation*}%
for every $u=\left( u_{1},\ldots ,u_{n}\right) $ and $v=\left( v_{1},\ldots
,v_{n}\right) \in \mathbb{R}^{n}$. When a metric is allowed to have index $1$%
, a trichotomy of vector classes named as \textbf{causality} arises. The
causality of $u\in \mathbb{\mathbb{L}}^{n}$ is said

\begin{description}
\item[i.] \textbf{spacelike}, when $\left\langle u,u\right\rangle _{1}>0$ or 
$u=0$;

\item[ii.] \textbf{timelike}, when $\left\langle u,u\right\rangle _{1}<0$;

\item[iii.] \textbf{lightlike}, when $u\neq 0$ and $\left\langle
u,u\right\rangle _{1}=0$.
\end{description}

The causality of an immersion $i:M\longrightarrow \mathbb{\mathbb{L}}^{n}$
is said \textbf{spacelike} (resp. \textbf{timelike} or \textbf{lightlike})
when all of its normal vectors are timelike (resp. spacelike or lightlike).
Finally, an immersion $i:M\longrightarrow \mathbb{\mathbb{L}}^{n}$ is called 
\textbf{nondegenarated}\textit{\ when it is} a spacelike or a timelike
immersion.

\bigskip

As usual, the \textbf{norm} of a vector $u\in \mathbb{L}^{n}$ is defined as $%
\left\Vert u\right\Vert _{1}=\sqrt{\left\vert \left\langle u,u\right\rangle
_{1}\right\vert }$. The Lorentzian \textbf{cross product} is the operator
determinated by

\begin{equation*}
u_{1}\times _{1}\ldots \times _{1}u_{n-1}=\det \left( 
\begin{array}{cccc}
e_{1} & \cdots & e_{n-1} & -e_{n} \\ 
u_{1}^{1} & \cdots & u_{1}^{n-1} & u_{1}^{n} \\ 
\vdots & \vdots & \vdots & \vdots \\ 
u_{n-1}^{1} & \cdots & u_{n-1}^{n-1} & u_{n-1}^{n}%
\end{array}%
\right) ,
\end{equation*}%
where $e_{1},\ldots ,e_{n}$ form the canonical basis of $\mathbb{R}^{n}$ and 
$u_{i}=\left( u_{i}^{1},\ldots ,u_{i}^{n}\right) $ are arbitrary vectors of $%
\mathbb{L}^{n}$, for every $i=1,\ldots ,n-1$.

\bigskip

A non degenerated smooth curve $\gamma :\left( a,b\right) \rightarrow 
\mathbb{L}^{3}$ is called a \textbf{Frenet curve} when\textbf{\ }it is a
biregular ($\left\Vert \gamma ^{\prime \prime }\right\Vert _{1}\neq 0$)
curve parametrized by arc lenght \TEXTsymbol{\backslash}($\left\Vert \gamma
^{\prime }\right\Vert _{1}\equiv 1$). In this case, we can define the 
\textbf{tangent vectors}, the \textbf{principal normal} vectors and \textbf{%
binormal vectors} of $\gamma $ by%
\begin{equation*}
T=\gamma ^{\prime }\;\text{,}\qquad N=\tfrac{\gamma ^{\prime \prime }}{%
\left\Vert \gamma ^{\prime \prime }\right\Vert _{1}}\qquad \text{and}\qquad
B=T\times _{1}N\text{.}
\end{equation*}%
The \textbf{curvature} and the \textbf{torsion} of a Frenet curve $\gamma $
are defined by%
\begin{equation*}
\kappa =\left\Vert \gamma ^{\prime \prime }\right\Vert _{1}\qquad \text{and}%
\qquad \tau =\left\langle N^{\prime },B\right\rangle _{1}
\end{equation*}%
and the Lorentzian Frenet formulas of $\gamma $ can be matricially expressed
as%
\begin{equation}
\left( 
\begin{array}{c}
T^{\prime } \\ 
N^{\prime } \\ 
B^{\prime }%
\end{array}%
\right) =\left( 
\begin{array}{ccc}
0 & \kappa  & 0 \\ 
-\varepsilon _{T}\varepsilon _{N}\kappa  & 0 & \tau  \\ 
0 & \varepsilon _{T}\tau  & 0%
\end{array}%
\right) \left( 
\begin{array}{c}
T \\ 
N \\ 
B%
\end{array}%
\right) ,  \label{frenet formula}
\end{equation}%
where $\varepsilon _{T}=\left\langle T,T\right\rangle _{1}$ and $\varepsilon
_{N}=\left\langle N,N\right\rangle _{1}$. Observe that the Lorentzian Frenet
formulas are natural generalizations for the well known Euclidean Frenet
formulas (because $\varepsilon _{T},\varepsilon _{N}\equiv 1$ in Euclidean $3
$-space).

Finally, observe that a Frenet curve $\gamma $ is spacelike (resp. timelike)
whenever $\varepsilon _{T}\equiv 1$ (resp. $\varepsilon _{T}\equiv -1$). In
this direction, we can define the \textbf{normal casuality} of $\gamma $ as 
\textbf{spacelike} (resp. \textbf{timelike}) whenever $\varepsilon
_{N}\equiv 1$ (resp. $\varepsilon _{N}\equiv -1$). For technical reasons, it
is also convenient to define $\varepsilon _{B}=\left\langle B,B\right\rangle
_{1}=-\varepsilon _{T}\varepsilon _{N}$. 

\bigskip

The \textbf{hyperbolic }$3$-\textbf{space} is the set

\begin{equation*}
\mathbb{H}^{3}=\left\{ x=\left( x_{1},\ldots ,x_{4}\right) \in \mathbb{L}^{4}%
\text{ };\text{ }\left\langle x,x\right\rangle _{1}=-1\text{ where }%
x_{4}>0\right\} 
\end{equation*}%
endowed with the metric induced by the canonical inclusion in $\mathbb{L}^{4}
$. Note that $\mathbb{H}^{3}$ is a Riemannian manifold since all tangent
vector of $\mathbb{H}^{3}$ are spacelike vectors.

\bigskip

Just like in the Lorentzian case, a smooth curve $\gamma :\left( a,b\right)
\rightarrow \mathbb{H}^{3}$ is called a \textbf{Frenet curve} when\textbf{\ }%
it is parametrized by arc lenght ($\left\Vert \gamma ^{\prime }\right\Vert
_{1}\equiv 1$) and verify $\left\Vert \gamma ^{\prime \prime }-\gamma
\right\Vert _{1}\neq 0$. In this case, we can also define the \textbf{%
tangent vector}, the \textbf{principal normal} vector and \textbf{binormal
vector} of $\gamma $ by%
\begin{equation*}
T=\gamma ^{\prime }\;\text{,}\qquad N=\tfrac{\gamma ^{\prime \prime }-\gamma 
}{\left\Vert \gamma ^{\prime \prime }-\gamma \right\Vert _{1}}\qquad \text{%
and}\qquad B=\gamma \times _{1}T\times _{1}N\text{.}
\end{equation*}%
The \textbf{curvature} and the \textbf{torsion} of a Frenet curve $\gamma $
are defined by%
\begin{equation*}
\kappa =\left\Vert \gamma ^{\prime \prime }-\gamma \right\Vert _{1}\qquad 
\text{and}\qquad \tau =\left\langle N^{\prime },B\right\rangle _{1}
\end{equation*}%
and the Hyperbolic Frenet formulas of $\gamma $ can be matricially expressed
as%
\begin{equation}
\left( 
\begin{array}{c}
\gamma ^{\prime } \\ 
T^{\prime } \\ 
N^{\prime } \\ 
B^{\prime }%
\end{array}%
\right) =\left( 
\begin{array}{cccc}
0 & 1 & 0 & 0 \\ 
1 & 0 & \kappa  & 0 \\ 
0 & -\kappa  & 0 & \tau  \\ 
0 & 0 & -\tau  & 0%
\end{array}%
\right) \left( 
\begin{array}{c}
\gamma  \\ 
T \\ 
N \\ 
B%
\end{array}%
\right)   \label{hyperbolic Franet formulas}
\end{equation}

\bigskip

Given a non degenerated regular surface $S$ in the Lorentzian $3$-space
(resp. a regular surface in the Hyperbolic $3$-space) locally parametrized
by an immersion $\psi \left( s,t\right) $, we define its canonical \textbf{%
normal vector field} of $S$ (associated to $\psi $) by%
\begin{equation*}
\mathcal{N}=\frac{\mathfrak{\psi }_{s}\times _{1}\mathfrak{\psi }_{t}}{%
\left\Vert \mathfrak{\psi }_{s}\times _{1}\mathfrak{\psi }_{t}\right\Vert
_{1}}\mathbb{\qquad }\left( \text{resp. }\mathcal{N}=\frac{\psi \times _{1}%
\mathfrak{\psi }_{s}\times _{1}\mathfrak{\psi }_{t}}{\left\Vert \psi \times
_{1}\mathfrak{\psi }_{s}\times _{1}\mathfrak{\psi }_{t}\right\Vert _{1}}%
\right) \text{.}
\end{equation*}%
The \textbf{Gaussian} and \textbf{mean curvatures} of $S$ are given by%
\begin{equation*}
K=\varepsilon \frac{eg-f^{2}}{EG-F^{2}}\qquad \text{and}\qquad H=\varepsilon 
\frac{eG-2fF+gE}{2\left( EG-F^{2}\right) }
\end{equation*}%
where 
\begin{eqnarray*}
E &=&\left\langle \mathfrak{\psi }_{s},\mathfrak{\psi }_{s}\right\rangle _{1}%
\text{,}\qquad F=\left\langle \mathfrak{\psi }_{s},\mathfrak{\psi }%
_{t}\right\rangle _{1}\text{,}\qquad G=\left\langle \mathfrak{\psi }_{t},%
\mathfrak{\psi }_{t}\right\rangle _{1}\text{,} \\
&& \\
e &=&\left\langle \mathfrak{\psi }_{ss},\mathcal{N}\right\rangle _{1}\text{,}%
\qquad f=\left\langle \mathfrak{\psi }_{st},\mathcal{N}\right\rangle _{1}%
\text{,}\qquad g=\left\langle \mathfrak{\psi }_{tt},\mathcal{N}\right\rangle
_{1}\text{,}
\end{eqnarray*}%
and%
\begin{equation*}
\varepsilon =\left\langle \mathcal{N},\mathcal{N}\right\rangle \in \left\{
-1,1\right\} .
\end{equation*}%
Observe that the number $\varepsilon $, called the \textbf{signal} of $S$,
is an attribute of the surface $S$, that is, it is independent of the
parametrization $\psi $. We remark that the curvatures presented above are
natural generalizations for the Euclidean ones (because $\varepsilon \equiv 1
$ in Euclidean space).

\section{Tubular Surfaces in $\mathbb{E}^{3}$, $\mathbb{L}^{3}$ and $\mathbb{%
H}^{3}$\label{Tubular Surface}}

\qquad Consider a smooth curve $\gamma :\left( a,b\right) \rightarrow 
\mathbb{E}^{3}$ (called \textbf{central curve}) with non zero tangent
vectors. A\textbf{\ tubular surface} of radius $r>0$ around $\gamma $ is the
set obtained by the union of all circles $S_{r}\left( \gamma \left( s\right)
\right) $ of radius $r$ and center $\gamma \left( s\right) $ contained in
the normal planes $T_{s}^{\bot }\gamma $ of $\gamma $.

A tubular surface is called \textbf{regular} whenever $\gamma $ is
parametrized by arc length and%
\begin{equation}
\xi \left( s,t\right) =1-r\kappa (s)\cos t\neq 0\text{,}
\label{regularidade}
\end{equation}%
for every $\left( s,t\right) \in \left( a,b\right) \times \mathbb{R}$. In
intervals where $\gamma $ is a Frenet curve, we observe that the technical
condition (\ref{regularidade}) is precisely the necessary and sufficiently
condition for the application%
\begin{equation}
\psi :\left( s,t\right) \in \left( a,b\right) \times \mathbb{R\longmapsto
\gamma }\left( s\right) +r\cos \left( t\right) N\left( s\right) +r\sin
\left( t\right) B\left( s\right) \in \mathbb{R}^{3}  \label{parametrization}
\end{equation}%
be an immersion. At such intervals, simple calculations with the above
parametrization provide the following expressions for the Gaussian, Mean and
principal curvatures of a regular tubular surface:

\begin{equation*}
K=\frac{\kappa \cos t}{r\left( r\kappa \cos t-1\right) }\,\text{,}\qquad H=%
\frac{2r\kappa \cos t-1}{2r\left( r\kappa \cos t-1\right) }\,\text{,}
\end{equation*}%
\begin{equation*}
k_{1}=\frac{\kappa \cos t}{r\kappa \cos t-1}\qquad \text{and}\qquad k_{2}=%
\frac{1}{r}\text{.}
\end{equation*}%
It is not difficult to see that the expressions above are also valid on
points where $\kappa (s)=0$.

\vspace{0.5cm}

In the Lorentzian space, tubular surfaces are defined in the same way, but
we also accept hyperbolic normal sections. More precisely, all the normal
sections can be Lorentzian circles%
\begin{equation*}
S_{r}\left( \gamma \left( s\right) \right) =\left\{ x\in T_{s}\gamma ^{\bot }%
\text{ };\text{ }\left\langle x-\gamma \left( s\right) ,x-\gamma \left(
s\right) \right\rangle _{1}=r^{2}\right\} \text{,}
\end{equation*}%
or all the normal sections can be Lorentzian hyperboles%
\begin{equation*}
S_{r}\left( \gamma \left( s\right) \right) =\left\{ x\in T_{s}\gamma ^{\bot }%
\text{ };\text{ }\left\langle x-\gamma \left( s\right) ,x-\gamma \left(
s\right) \right\rangle _{1}=-r^{2}\right\} \text{.}
\end{equation*}%
Note that Lorentzian hyperboles are Euclidean hyperboles, however Lorentzian
circles can not be Euclidean circles.

\bigskip

In intervals where the central curve $\gamma :\left( a,b\right) \rightarrow 
\mathbb{L}^{3}$ is a Frenet curve, the tubular surface admits a global
parameterization that depends on the causality and the normal causality of $%
\gamma $. The table below contain all possibilities:%
\begin{equation}
\begin{tabular}{||c||c||c||c||}
\hline\hline
Curve & Normal & Section & Parametrization \\ \hline\hline
spacelike & spacelike & Lorentzian circles & $\gamma \pm r\cosh \left(
t\right) N+r\sinh \left( t\right) B$ \\ \hline\hline
spacelike & timelike & Lorentzian circles & $\gamma +r\sinh \left( t\right)
N\pm r\cosh \left( t\right) B$ \\ \hline\hline
timelike & spacelike & Lorentzian circles & $\gamma +r\cos \left( t\right)
N+r\sin \left( t\right) B$ \\ \hline\hline
spacelike & spacelike & Lorentzian hyperboles & $\gamma +r\sinh \left(
t\right) N\pm r\cosh \left( t\right) B$ \\ \hline\hline
spacelike & timelike & Lorentzian hyperboles & $\gamma \pm r\cosh \left(
t\right) N+r\sinh \left( t\right) B$ \\ \hline\hline
timelike & spacelike & Lorentzian hyperboles & does not exist \\ \hline\hline
\end{tabular}
\label{table of parametrizations}
\end{equation}

\bigskip

In the aim to encompass and deal with all possible local parametrizations
simultaneosly, we will use a pair of functions $\left( \mu \left( t\right)
,\eta \left( t\right) \right) $ to represent one of the followings pairs%
\begin{equation}
\left( \delta \cos t,\sin t\right) \text{, }\left( \delta \cosh t,\sinh
t\right) \text{,}\ \left( \sinh t,\delta \cosh t\right) \text{,}
\label{mi eta cjto 3}
\end{equation}%
where $\delta \in \left\{ -1,1\right\} $. Then, each tubular surface
presented in (\ref{table of parametrizations}) is rewritten as%
\begin{equation}
\psi \left( s,t\right) =\gamma \left( s\right) +r\mu \left( t\right) N\left(
s\right) +r\eta \left( t\right) B\left( s\right) \text{,}
\label{generic param L 1}
\end{equation}%
for a given pair $\left( \mu \left( t\right) ,\eta \left( t\right) \right) $
in (\ref{mi eta cjto 3}). Observe that 
\begin{equation}
\mu ^{\prime }\left( t\right) =\delta \varepsilon _{T}\eta \left( t\right)
\qquad \text{and}\qquad \eta ^{\prime }\left( t\right) =\delta \mu \left(
t\right) \text{.}  \label{hip1}
\end{equation}

\bigskip

Generalizing the Euclidean case, we say that a tubular surface in the
Lorentzian space is \textbf{regular} when $\gamma $ is parametrized by arc
length and%
\begin{equation}
\xi \left( s,t\right) =1+\varepsilon _{B}r\kappa \left( s\right) \mu \left(
t\right) \neq 0\text{,}  \label{regularidade_Lorentziano}
\end{equation}%
for every $\left( s,t\right) \in \left( a,b\right) \times \mathbb{R}$. In
intervals where the central curve $\gamma $ is a Frenet curve, it follows
from (\ref{frenet formula}) and (\ref{generic param L 1}) that%
\begin{equation*}
\psi _{t}\left( s,t\right) =r\varepsilon _{T}\eta \left( t\right) N\left(
s\right) +r\mu \left( t\right) B\left( s\right)
\end{equation*}%
and%
\begin{equation*}
\psi _{s}\left( s,t\right) =\xi \left( s,t\right) T\left( s\right) +\delta
\tau \psi _{t}\left( s,t\right) \text{,}
\end{equation*}%
for every $\left( s,t\right) \in \left( a,b\right) \times \mathbb{R}$. Then,
it is easy to see that $\psi $ is an immersion if and only if condition (\ref%
{regularidade_Lorentziano}) holds. Finally, at such intervals, an unit
normal vector of the tubular surface is given by%
\begin{equation*}
\mathcal{N}\left( s,t\right) =-\mu \left( t\right) N\left( s\right) -\eta
\left( t\right) B\left( s\right) \text{.}
\end{equation*}%
Standard calculations with the above parametrization provide the following
expressions for the Gaussian and Mean curvatures of a regular tubular
surface in the Lorentzian space:

\begin{equation}
K=\varepsilon \dfrac{\varepsilon _{B}\kappa \left( s\right) \mu \left(
t\right) }{r\left( 1+\varepsilon _{B}r\kappa \left( s\right) \mu \left(
t\right) \right) }\mathbb{\qquad }\text{and}\qquad H=\varepsilon \dfrac{%
2\varepsilon _{B}r\kappa \left( s\right) \mu \left( t\right) +1}{2r\left(
1+\varepsilon _{B}r\kappa \left( s\right) \mu \left( t\right) \right) }\text{%
,}  \label{Tubular_curvatures}
\end{equation}%
where $\varepsilon $ is signal of the surface parametrized by $\psi $.
Again, the expressions above are also valid in points where $\kappa \left(
s\right) =0$.

\bigskip

\qquad In the Hyperbolic $3$-space,\textbf{\ }tubular surfaces\textbf{\ }%
have the same definition as in the Euclidean case. In intervals where $%
\gamma $ is a Frenet curve, it is not difficult show that tubular surfaces
with Frenet central curves admit a parametrization of the form%
\begin{equation}
\psi \left( s,t\right) =\left( \cosh r\right) \gamma +\sinh r\left( \cos
tN+\sin tB\right) \text{.}  \label{parametrization hyperbolic}
\end{equation}%
In this case, the regularity condition of the tubular surfaces is%
\begin{equation*}
\xi _{\mathbb{H}}\left( s,t\right) =\cosh r-\kappa \cos t\sinh r\neq 0\text{,%
}
\end{equation*}%
for every $\left( s,t\right) \in \left( a,b\right) \times \mathbb{R}$. We
have the following expressions for the Gaussian, Mean and principal
curvatures of a regular tubular surface:

\begin{equation*}
K=-\frac{\kappa \cos t}{\xi _{\mathbb{H}}\sinh r}\;,\mathbb{\qquad }H=\frac{%
\cosh r-2\kappa \cos t\sinh r}{2\xi _{\mathbb{H}}\sinh r}\;
\end{equation*}%
\begin{equation*}
k_{1}=-\frac{\kappa \cos t}{\xi _{\mathbb{H}}}\qquad \text{and}\qquad k_{2}=%
\frac{1}{\sinh r}\text{.}
\end{equation*}

Finally, we say that a (non degenerated) regular tubular surface in any of
the three spaces mentioned above is a \textbf{right cilynder} when the
central curve has zero curvature ($\kappa \equiv 0$), that is, when the
central curve is a geodesic of the ambient.

\section{Polynomials results for Tubular Surfaces\label{Polynomials results}}

The objective of this technical session is to obtain (a reasonably simple to
check) condition to determine when a polynomial $Q\left( x,y\right) \in 
\mathbb{R}\left[ x,y\right] $ has a factor of the form $xr^{2}-2yr+1$, where 
$r\in \mathbb{R}-\left\{ 0\right\} $. The interest in studying this specific
type of polynomial naturally arises in the investigation of polynomial
Weingarten tubular surfaces (of radius $r$) and will become clearer in the
next section. As a matter of fact, we obtain the following theorem.

\begin{theorem}
\label{Teo_Factorization}Consider a polynomial $Q\left( x,y\right) \in 
\mathbb{R}\left[ x,y\right] $ and let $r$ be nonzero constant. Then $Q$
belongs to the ideal generated by $xr^{2}-2yr+1$ if, and only if, $Q\left( 
\frac{x}{r},\frac{xr+1}{2r}\right) \in \mathbb{R}\left[ x\right] $ is
identically null.
\end{theorem}

\bigskip

The following numerical example will illustrate one first manner to apply
the theorem above.

\begin{example}
\label{Ex Q}Consider the polynomial%
\begin{eqnarray*}
Q\left( x,y\right) 
&=&4x^{4}+8x^{2}y^{2}-12xy^{3}+9x^{3}+9x^{2}y-9xy^{2}-4y^{3} \\
&&\mathbb{\qquad }+22x^{2}-8xy-7y^{2}-91x+98y-24
\end{eqnarray*}%
and notice that 
\begin{eqnarray*}
Q\left( \frac{x}{r},\frac{xr+1}{2r}\right)  &=&\left( \frac{-3r^{3}+4r^{2}+8%
}{2r^{4}}\right) x^{4}-\left( \frac{2r^{3}+9r^{2}-52}{4r^{3}}\right) x^{3} \\
&&-\left( \frac{7r^{4}+22r^{3}-70r^{2}-8}{4r^{4}}\right) x^{2} \\
&&-\left( \frac{-196r^{4}+378r^{3}+22r^{2}+9r+6}{4r^{4}}\right) x \\
&&-\left( \frac{96r^{3}-196r^{2}+7r+2}{4r^{3}}\right) .
\end{eqnarray*}%
$\allowbreak \allowbreak $Then, for $r=2$, we have%
\begin{equation*}
Q\left( \frac{x}{2},\frac{2x+1}{4}\right) \equiv 0.
\end{equation*}%
By Theorem \ref{Teo_Factorization}, the polynomial $Q\left( x,y\right) $ can
be divided by $4x-4y+1$. In fact, polynomial $Q\left( x,y\right) $ can be
rewritten as%
\begin{equation*}
Q\left( x,y\right) =\left( 4x-4y+1\right) \left(
x^{3}+x^{2}y+3xy^{2}+2x^{2}+4xy+y^{2}+5x+2y-24\right) .
\end{equation*}
\end{example}

\bigskip

The first difficulty we face when we work with polynomials is to find a
favorable expression that allows us to examine them. In this direction, we
present a Propostion that provides an explicit expression for the
coefficients of the polynomial $Q\left( \frac{x}{r},\frac{xr+1}{2r}\right)
\in \mathbb{R}\left[ x\right] $. In the proof we will use the following
(easy to check) identities:%
\begin{equation}
\sum\limits_{i=0}^{n+1}\sum\limits_{j=0}^{n+1-i}M_{i,j}=\sum%
\limits_{i=0}^{n+1}M_{i,n+1-i}+\sum\limits_{i=0}^{n}\sum%
\limits_{j=0}^{n-i}M_{i,j}  \label{Sumatorio_1}
\end{equation}%
\begin{equation}
\sum\limits_{k=0}^{n+1}\sum\limits_{i=0}^{k}\sum%
\limits_{j=0}^{n+1-k}M_{i,j,k}=\sum\limits_{k=0}^{n}\sum\limits_{i=0}^{k}%
\sum\limits_{j=0}^{n-k}M_{i,j,k}+\sum\limits_{k=0}^{n+1}\sum%
\limits_{i=0}^{k}M_{i,n+1-k,k}  \label{Sumatorio_5}
\end{equation}%
\begin{equation}
\sum\limits_{k=0}^{n+1}\sum\limits_{i=0}^{k}M_{k,i}=\sum\limits_{i=0}^{n+1}%
\sum\limits_{k=0}^{n+1-i}M_{k+i,i}  \label{Sumatorio_4}
\end{equation}

\begin{proposition}
\label{Prop Polynomial 1}Consider a polynomial $Q\left( x,y\right)
=\sum\limits_{i=0}^{n}\sum\limits_{j=0}^{n-i}a_{i,j}x^{i}y^{j}\in \mathbb{R}%
\left[ x,y\right] $. For every $r\in \mathbb{R}-\left\{ 0\right\} $, we have 
\begin{equation*}
Q\left( \frac{x}{r},\frac{xr+1}{2r}\right) =\sum\limits_{k=0}^{n}\Gamma
_{k}^{n}\left( r\right) x^{k}\text{,}
\end{equation*}%
where%
\begin{equation}
\Gamma _{k}^{n}\left( r\right) =\sum\limits_{i=0}^{k}\sum\limits_{j=0}^{n-k}%
\tbinom{k-i+j}{j}\tfrac{a_{i,k-i+j}}{2^{k-i+j}r^{j+i}}\text{.}
\label{def Tau}
\end{equation}
\end{proposition}

\begin{proof}
Going along the lines of the proof, we will demonstrate%
\begin{equation}
\Psi _{n}\left( r\right)
:=\sum\limits_{i=0}^{n}\sum\limits_{j=0}^{n-i}a_{i,j}\left( \tfrac{x}{r}%
\right) ^{i}\left( \tfrac{xr+1}{2r}\right) ^{j}-\sum\limits_{k=0}^{n}\Gamma
_{k}^{n}\left( r\right) x^{k}=0\text{,}  \label{pol expression}
\end{equation}%
by induction on $n$. It is easy to see that the above equality is verified
for $n\in \left\{ 0,1\right\} $.

\bigskip

So assume that the relation (\ref{pol expression}) is true for some $n\in 
\mathbb{N}$. Replacing $\Gamma _{k}^{n+1}\left( r\right) $ by the expression
in (\ref{def Tau}) and using identity \ref{Sumatorio_1} we obtain%
\begin{eqnarray*}
\Psi _{n+1}\left( r\right)
&=&\sum\limits_{i=0}^{n+1}\sum\limits_{j=0}^{n+1-i}a_{i,j}\left( \tfrac{x}{r}%
\right) ^{i}\left( \tfrac{xr+1}{2r}\right)
^{j}-\sum\limits_{k=0}^{n+1}\Gamma _{k}^{n+1}\left( r\right) x^{k} \\
&=&\sum\limits_{i=0}^{n+1}a_{i,n+1-i}\left( \tfrac{x}{r}\right) ^{i}\left( 
\tfrac{xr+1}{2r}\right)
^{n+1-i}+\sum\limits_{i=0}^{n}\sum\limits_{j=0}^{n-i}a_{i,j}\left( \tfrac{x}{%
r}\right) ^{i}\left( \tfrac{xr+1}{2r}\right) ^{j} \\
&&\mathbb{\qquad }\qquad
-\sum\limits_{k=0}^{n+1}\sum\limits_{i=0}^{k}\sum\limits_{j=0}^{n+1-k}%
\tbinom{k-i+j}{j}\tfrac{a_{i,k-i+j}}{2^{k-i+j}r^{j+i}}x^{k}.
\end{eqnarray*}%
So the identity \ref{Sumatorio_5} applied to the last term and the induction
hypothesis give us%
\begin{equation*}
\Psi _{n+1}\left( r\right) =\sum\limits_{i=0}^{n+1}a_{i,n+1-i}\left( \tfrac{x%
}{r}\right) ^{i}\left( \tfrac{xr+1}{2r}\right)
^{n+1-i}-\sum\limits_{k=0}^{n+1}\sum\limits_{i=0}^{k}\tbinom{n+1-i}{n+1-k}%
\tfrac{a_{i,n+1-i}}{2^{n+1-i}r^{n+1+i-k}}x^{k}
\end{equation*}%
Finally, applying the Binomial Theorem to the first term and identity \ref%
{Sumatorio_4} to the last one, we have%
\begin{eqnarray*}
\Psi _{n+1}\left( r\right)
&=&\sum\limits_{i=0}^{n+1}\sum\limits_{k=0}^{n+1-i}\tbinom{n+1-i}{k}\tfrac{%
a_{i,n+1-i}}{2^{n+1-i}r^{n+1-k}}x^{k+i} \\
&&\mathbb{\qquad }\qquad -\sum\limits_{i=0}^{n+1}\sum\limits_{k=0}^{n+1-i}%
\tbinom{n+1-i}{n+1-i-k}\tfrac{a_{i,n+1-i}}{2^{n+1-i}r^{n+1-k}}x^{k+i} \\
&=&0.
\end{eqnarray*}%
\hfill
\end{proof}

\bigskip

Next example present how Proposition \ref{Prop Polynomial 1} allows us to
gather more information about a polynomial $Q$ and explain how to choose $r$
to test whether $Q\left( \frac{x}{r},\frac{xr+1}{2r}\right) $ is identically
null or not.

\begin{example}
\label{Exemplo2}Consider again the polynomial $Q\left( x,y\right) $
presented in the Example \ref{Ex Q}. By Proposition \ref{Prop Polynomial 1}
the independent term of $Q\left( \tfrac{x}{r},\tfrac{xr+1}{2r}\right) \ $is
given by%
\begin{eqnarray*}
\Gamma _{0}^{4}\left( r\right) &=&\frac{98}{2r}+\frac{-7}{4r^{2}}+\frac{-4}{%
8r^{3}}+\frac{0}{16r^{4}}-24 \\
&=&-\frac{24\left( r-\frac{1}{8}\right) \left( r-2\right) \left( r+\frac{1}{%
12}\right) }{r^{3}}.
\end{eqnarray*}%
Since $\Gamma _{0}^{4}\left( 1\right) \neq 0$, Theorem \ref%
{Teo_Factorization} allow us to say that $Q\left( x,\frac{x+1}{2}\right) $
does not belong to the ideal generated by $x-2y+1$.
\end{example}

\bigskip

The next proposition is a technical result that presents us a
characterization of polynomials $Q\left( x,y\right) $ that belongs to the
ideal generated by $xr^{2}-2yr+1$.

\begin{proposition}
\label{Prop Polynomial 2}Consider a real number $r\neq 0$. A polynomial $%
Q\left( x,y\right)
=\sum\limits_{i=0}^{n}\sum\limits_{j=0}^{n-i}a_{i,j}x^{i}y^{j}\in \mathbb{R}%
\left[ x,y\right] $ is divisible by $xr^{2}-2yr+1$ if, and only if, there
exists a family%
\begin{equation*}
\mathcal{C}=\left\{ c_{i,j}\in \mathbb{R}\;;\;i,j\in \mathbb{\mathbb{Z}}%
\right\}
\end{equation*}%
such that:

\begin{enumerate}
\item $c_{i,j}=0$ whenever $i<0$ or $j<0$ or $i+j\geq n$;

\item $a_{i,j}=r^{2}c_{i-1,j}-2rc_{i,j-1}+c_{i,j}$, for every $i,j\in 
\mathbb{\mathbb{Z}}.$
\end{enumerate}
\end{proposition}

\begin{proof}
Suppose that $Q\left( x,y\right) $ is divisible by $xr^{2}-2yr+1$. Then,
there is a family $\left\{ c_{i,j}\in \mathbb{R}\;;\;i\text{, }j\in \mathbb{N%
}\text{ and }i+j\leq n-1\right\} $ such that%
\begin{equation}
Q\left( x,y\right) =\left( xr^{2}-2yr+1\right)
\sum\limits_{i=0}^{n-1}\sum\limits_{j=0}^{n-1-i}c_{i,j}X^{i}Y^{j}\text{.}
\label{Q div}
\end{equation}%
It follows that%
\begin{equation*}
a_{i,j}=r^{2}c_{i-1,j}-2rc_{i,j-1}+c_{i,j}
\end{equation*}%
for every $i=1,\ldots ,n$ and for every $j=1,\ldots ,n-i$. If we set $%
c_{i,j}=0$, whenever $i<0$ or $j<0$ or $i+j\geq n$, we constructed a family $%
\mathcal{C}$ verifying items $1$ and $2$.

\bigskip

Conversely, given a family $\mathcal{C}$ as in the statement of the
proposition, define the polynomial 
\begin{equation*}
R\left( x,y\right)
=\sum\limits_{i=0}^{n-1}\sum\limits_{j=0}^{n-1-i}c_{i,j}x^{i}y^{j}.
\end{equation*}
It is easy to check that $Q\left( x,y\right) =\left( xr^{2}-2yr+1\right)
R\left( x,y\right) $.
\end{proof}

\bigskip

Before stating the main result of the section, we present a (probably not
original) lemma with a binomial identity that will be necessary to its
proof. Remind that (Pascal's Formula) 
\begin{equation*}
\dbinom{n+1}{k+1}=\dbinom{n}{k+1}+\dbinom{n}{k}\text{,}
\end{equation*}%
for every $n,k\in \mathbb{\mathbb{N}}$.

\begin{lemma}
\label{Prop 10.20}Given $n\in \mathbb{N}$ and $x,j\in \mathbb{\mathbb{Z}}$,
we have%
\begin{equation}
\sum_{m=0}^{n}\left( -1\right) ^{m}\dbinom{x-m}{j}\binom{n}{m}=\dbinom{x-n}{%
j-n}\text{.}  \label{Id_Lemma}
\end{equation}
\end{lemma}

\begin{proof}
For $n=1$, expression (\ref{Id_Lemma}) is a straightforward consequence of
Pascal's Formula. Suppose expression (\ref{Id_Lemma}) is true for some $n\in 
\mathbb{N}$. Then, Pascal's Formula implies 
\begin{equation*}
\sum_{m=0}^{n+1}\left( -1\right) ^{m}\tbinom{x-m}{j}\tbinom{n+1}{m}%
=\sum_{m=0}^{n}\left( -1\right) ^{m}\tbinom{x-m}{j}\tbinom{n}{m}%
-\sum_{m=0}^{n}\left( -1\right) ^{m}\tbinom{x-1-m}{j}\tbinom{n}{m}
\end{equation*}%
Applying the induction step for both terms on the right hand side yields%
\begin{equation*}
\sum_{m=0}^{n+1}\left( -1\right) ^{m}\tbinom{x-m}{j}\tbinom{n+1}{m}=\tbinom{%
x-n}{j-n}-\tbinom{x-n-1}{j-n}\text{.}
\end{equation*}%
Note that we can not apply Pascal's Formula directly on the above expression
to conclude that%
\begin{equation}
\tbinom{x-n}{j-n}-\tbinom{x-n-1}{j-n}=\tbinom{x-\left( n+1\right) }{j-\left(
n+1\right) }\text{.}  \label{Bin_identity}
\end{equation}%
There are some cases to consider:

\bigskip

When $j-n\leq 0$, identity (\ref{Bin_identity}) is true by definition. Then,
we can suppose that $j-n\geq 1$. To conclude (\ref{Bin_identity}) by
Pascal's Formula, we must have $x-n-1\geq 0$, that is, we must have $x\geq j$%
. Therefore, we only need to treat the case where $x<j$. But this case is
immediate beacause $x<j$ implies%
\begin{equation*}
\tbinom{x-n}{j-n}=\tbinom{x-n-1}{j-n}=\tbinom{x-\left( n+1\right) }{j-\left(
n+1\right) }=0
\end{equation*}%
by definition. Then, identity (\ref{Bin_identity}) is true in all cases and
this fact completes the proof.\hfill
\end{proof}

\bigskip

Theorem \ref{Teo_Factorization} is an immediate consequence of the following
theorem that articulates both of the previous propositions. In the proof we
will use the following (easy to check) identities:

\begin{equation}
\sum\limits_{m=0}^{i}\sum\limits_{l=0}^{n-i+m}\sum%
\limits_{k=0}^{i-m}M_{k,m,l}=\sum\limits_{k=0}^{i}\sum\limits_{m=0}^{i-k}%
\sum\limits_{l=0}^{n-i+m}M_{k,m,l}\text{,}i\in \left\{ 0,...,n\right\}
\label{E1}
\end{equation}

\begin{equation}
\sum\limits_{i=0}^{n}M_{i}=\sum\limits_{i=0}^{n}M_{n-i}  \label{inv_cont_0}
\end{equation}

\begin{equation}
\sum\limits_{m=0}^{p}\sum\limits_{l=0}^{q+p-m}\lambda
_{m,l}A_{m+l}=\sum\limits_{l=0}^{p}\sum\limits_{m=0}^{l}\lambda
_{m,l-m}A_{l}+\sum\limits_{l=1}^{q}\sum\limits_{m=0}^{p}\lambda
_{m,p+l-m}A_{p+l}\;,  \label{insane}
\end{equation}%
for every $p$ and $q,n\in \mathbb{N}$.

\begin{theorem}
\label{Prop Polynomial 3}Given $n\in \mathbb{N-}\left\{ 0\right\} $,
consider a family%
\begin{equation*}
\mathcal{A}=\left\{ a_{i,j}\in \mathbb{R}\;;\;i,j\in \mathbb{N}\right\}
\end{equation*}%
such that $a_{i,j}=0$, whenever either $i<0$ or $j<0$ or $i+j>n$. Then the
coefficients $a_{i,j}$ verify%
\begin{equation}
\Gamma _{k}^{n}\left( r\right) =\sum\limits_{i=0}^{k}\sum\limits_{j=0}^{n-k}%
\binom{k-i+j}{j}\frac{a_{i,k-i+j}}{2^{k-i+j}r^{j+i}}=0\qquad ,\forall
k=0,\ldots ,n  \label{Def Tau 2}
\end{equation}%
if and only if there is a family%
\begin{equation*}
\mathcal{C}=\left\{ c_{i,j}\in \mathbb{R}\;;\;i,j\in \mathbb{Z}\right\}
\end{equation*}%
verifying

\begin{enumerate}
\item $c_{i,j}=0$ whenever $i<0$ or $j<0$ or $i+j\geq n$;

\item $a_{i,j}=r^{2}c_{i-1,j}-2rc_{i,j-1}+c_{i,j}$, for every $i,j\in 
\mathbb{N}$
\end{enumerate}
\end{theorem}

\begin{proof}
Assume the existence of a family $\mathcal{C}\ $as stated in the theorem. We
need to prove that $\Gamma _{k}^{n}\left( r\right) =0$, for every $k$.

\bigskip

When $k\in \left\{ 0,n\right\} $, the desired assertion follows from
properties (1) and (2). In fact:%
\begin{eqnarray*}
\Gamma _{0}^{n}\left( r\right) &=&\sum\limits_{j=0}^{n}\tfrac{a_{0,j}}{%
2^{j}r^{j}}=\sum\limits_{j=0}^{n}\tfrac{-2rc_{0,j-1}+c_{0,j}}{2^{j}r^{j}} \\
&=&\sum\limits_{j=0}^{n-1}\tfrac{c_{0,j}}{2^{j}r^{j}}-\sum\limits_{j=1}^{n}%
\tfrac{c_{0,j-1}}{2^{j-1}r^{j-1}}=0
\end{eqnarray*}%
\begin{eqnarray*}
\Gamma _{n}^{n}\left( r\right) &=&\sum\limits_{i=0}^{n}\tfrac{a_{i,n-i}}{%
2^{n-i}r^{i}}=\sum\limits_{i=0}^{n}\tfrac{r^{2}c_{i-1,n-i}-2rc_{i,n-i-1}}{%
2^{n-i}r^{i}} \\
&=&\sum\limits_{i=1}^{n}\tfrac{c_{i-1,n-i}}{2^{n-i}r^{i-2}}%
-\sum\limits_{i=0}^{n-1}\tfrac{c_{i,n-i-1}}{2^{n-i-1}r^{i-1}}=0
\end{eqnarray*}

\bigskip

Given $k\in \left\{ 1,\ldots ,n-1\right\} $, it follows from property (2)
that%
\begin{equation*}
\Gamma _{k}^{n}\left( r\right) =\sum\limits_{i=0}^{k}\sum\limits_{j=0}^{n-k}%
\tbinom{k-i+j}{j}\tfrac{r^{2}c_{i-1,k-i+j}-2rc_{i,k-i+j-1}+c_{i,k-i+j}}{%
2^{k-i+j}r^{j+i}}\text{.}
\end{equation*}%
Since $k\geq 1$, property (1) give us%
\begin{eqnarray*}
\Gamma _{k}^{n}\left( r\right)
&=&\sum\limits_{i=1}^{k}\sum\limits_{j=0}^{n-k}\tbinom{k-i+j}{j}\tfrac{%
c_{i-1,k-i+j}}{2^{k-i+j}r^{j+i-2}}-\sum\limits_{i=0}^{k-1}\sum%
\limits_{j=0}^{n-k}\tbinom{k-i+j}{j}\tfrac{c_{i,k-i+j-1}}{%
2^{k-i+j-1}r^{j+i-1}} \\
&&\mathbb{\qquad }-\sum\limits_{j=0}^{n-k}\tfrac{c_{k,j-1}}{2^{j-1}r^{j+k-1}}%
+\sum\limits_{i=0}^{k-1}\sum\limits_{j=0}^{n-k}\tbinom{k-i+j}{j}\tfrac{%
c_{i,k-i+j}}{2^{k-i+j}r^{j+i}}+\sum\limits_{j=0}^{n-k}\tfrac{c_{k,j}}{%
2^{j}r^{j+k}} \\
&=&\sum\limits_{i=0}^{k-1}\sum\limits_{j=0}^{n-k}\tbinom{k-i+j-1}{j}\tfrac{%
c_{i,k-i+j-1}}{2^{k-i+j-1}r^{j+i-1}}-\sum\limits_{i=0}^{k-1}\sum%
\limits_{j=0}^{n-k}\tbinom{k-i+j}{j}\tfrac{c_{i,k-i+j-1}}{%
2^{k-i+j-1}r^{j+i-1}} \\
&&\mathbb{\qquad }+\sum\limits_{i=0}^{k-1}\sum\limits_{j=0}^{n-k}\tbinom{%
k-i+j}{j}\tfrac{c_{i,k-i+j}}{2^{k-i+j}r^{j+i}} \\
&=&\sum\limits_{i=0}^{k-1}\sum\limits_{j=1}^{n-k}\left[ \tbinom{k-i+j-1}{j}-%
\tbinom{k-i+j}{j}+\tbinom{k-i+j-1}{j-1}\right] \tfrac{c_{i,k-i+j-1}}{%
2^{k-i+j-1}r^{j+i-1}}.
\end{eqnarray*}%
Then, $\Gamma _{k}^{n}\left( r\right) =0$ from Pascal's Formula.

\bigskip

Now suppose that $\Gamma _{k}^{n}\left( r\right) =0$, for every $k=\left\{
0,\ldots ,n\right\} $, and define the family $\mathcal{C}=\left\{ c_{i,j}\in 
\mathbb{R}\;;\;i,j\in \mathbb{Z}\right\} $ by%
\begin{equation}
c_{i,j}=\left\{ 
\begin{array}{ll}
\sum\limits_{k=0}^{i}\sum\limits_{l=0}^{j}\left( -1\right) ^{k}\binom{l+k}{l}%
2^{l}r^{l+2k}a_{i-k,j-l} & \text{,}\mathbb{\qquad }\text{if }i,j\in \mathbb{N%
} \\ 
0 & \text{,}\mathbb{\qquad }\text{if }i<0\text{ or }j<0%
\end{array}%
\right. \text{.}  \label{def cij}
\end{equation}%
With the aim of demonstrating that family $\mathcal{C}$ verify property (2)
define%
\begin{equation*}
\Phi _{i,j}\left( r\right) =r^{2}c_{i-1,j}-2rc_{i,j-1}+c_{i,j},
\end{equation*}%
for every $i,j\in \mathbb{Z}$. We have four cases to consider:

\bigskip

It is easy to see that $\Phi _{0,0}\left( r\right) =a_{0,0}$. For every $%
i\geq 1$, we have%
\begin{eqnarray*}
\Phi _{i,0}\left( r\right) &=&r^{2}\sum\limits_{k=0}^{i-1}\left( -1\right)
^{k}r^{2k}a_{i-1-k,0}+\sum\limits_{k=0}^{i}\left( -1\right)
^{k}r^{2k}a_{i-k,0} \\
&=&-\sum\limits_{k=1}^{i}\left( -1\right)
^{k}r^{2k}a_{i-k,0}+\sum\limits_{k=0}^{i}\left( -1\right) ^{k}r^{2k}a_{i-k,0}
\\
&=&a_{i,0}
\end{eqnarray*}%
Analogously, for every $j\geq 1$ we have $\Phi _{0,j}\left( r\right)
=a_{0,j} $. Finally, take $i,j\geq 1$ and observe that%
\begin{eqnarray*}
\Phi _{i,j}\left( r\right)
&=&r^{2}\sum\limits_{k=0}^{i-1}\sum\limits_{l=0}^{j}\left( -1\right) ^{k}%
\tbinom{l+k}{l}2^{l}r^{l+2k}a_{\left( i-1\right) -k,j-l} \\
&&\mathbb{\qquad }-2r\sum\limits_{k=0}^{i}\sum\limits_{l=0}^{j-1}\left(
-1\right) ^{k}\tbinom{l+k}{l}2^{l}r^{l+2k}a_{i-k,\left( j-1\right) -l} \\
&&\mathbb{\qquad }+\sum\limits_{k=0}^{i}\sum\limits_{l=0}^{j}\left(
-1\right) ^{k}\tbinom{l+k}{l}2^{l}r^{l+2k}a_{i-k,j-l} \\
&=&-\sum\limits_{k=1}^{i}\left( -1\right)
^{k}r^{2k}a_{i-k,j}+\sum\limits_{k=1}^{i}\sum\limits_{l=1}^{j}\left(
-1\right) ^{\left( k-1\right) }\tbinom{l+\left( k-1\right) }{l}%
2^{l}r^{l+2k}a_{i-k,j-l} \\
&&\mathbb{\qquad }-\sum\limits_{l=1}^{j}2^{l}r^{l}a_{i,j-l}-\sum%
\limits_{k=1}^{i}\sum\limits_{l=1}^{j}\left( -1\right) ^{k}\tbinom{\left(
l-1\right) +k}{\left( l-1\right) }2^{l}r^{l+2k}a_{i-k,j-l}+a_{i,j} \\
&&\mathbb{\qquad }+\sum\limits_{l=1}^{j}2^{l}r^{l}a_{i,j-l}+\sum%
\limits_{k=1}^{i}\left( -1\right) ^{k}r^{2k}a_{i-k,j} \\
&&\mathbb{\qquad }+\sum\limits_{k=1}^{i}\sum\limits_{l=1}^{j}\left(
-1\right) ^{k}\tbinom{l+k}{l}2^{l}r^{l+2k}a_{i-k,j-l} \\
&=&a_{i,j}+\sum\limits_{k=1}^{i}\sum\limits_{l=1}^{j}\left( -1\right)
^{\left( k-1\right) }\left[ \tbinom{l+\left( k-1\right) }{l}-\tbinom{\left(
l-1\right) +k}{\left( l-1\right) }+\tbinom{l+k}{l}\right]
2^{l}r^{l+2k}a_{i-k,j-l}.
\end{eqnarray*}%
Follows from Pascal's Formula that $\Phi _{i,j}\left( r\right) =a_{i,j}$.
Property (2) is therefore valid.

\bigskip

Now let's turn our attention to property (1). By definition (\ref{def cij}),
it remains to prove that $c_{i,j}=0$, for every $i,j\in \mathbb{N}$ such
that $i+j\geq n$. Let us split the proof into two claims:

\bigskip

\noindent \textbf{Claim 1:} $c_{i,j}=0$, for every $i,j\in \mathbb{N}$ such
that $i+j=n$.

\bigskip

In the presence of the hypothesis (\ref{Def Tau 2}), it is sufficient to
show that each $c_{i,n-i}$ is a linear combination of the coeficients $%
\Gamma _{k}^{n}\left( r\right) $. More precisely, it is sufficient to show
that%
\begin{equation*}
c_{i,n-i}=\sum\limits_{m=0}^{i}\left( -1\right) ^{m}\tbinom{n-i+m}{n-i}%
2^{n-i}r^{n+m}\Gamma _{i-m}^{n}\left( r\right)
\end{equation*}%
for every $i\in \left\{ 1,\ldots ,n\right\} $. Using the definition of the
coeficients $\Gamma _{k}^{n}\left( r\right) $ and the identity (\ref{E1}) we
conclude that is sufficient to show that%
\begin{equation}
c_{i,n-i}=\sum\limits_{k=0}^{i}\sum\limits_{m=0}^{i-k}\sum%
\limits_{l=0}^{n-i+m}\left( -1\right) ^{m}\frac{2^{k+m+n}r^{m+n}}{%
2^{l+2i}r^{l+k}}\tbinom{n-i+m}{n-i}\tbinom{i-m-k+l}{l}a_{k,i-m-k+l}
\label{eqpaiva}
\end{equation}%
for every $i\in \left\{ 1,\ldots ,n\right\} $. Applying identity (\ref%
{inv_cont_0}) to the definition (\ref{def cij}) we obtain%
\begin{equation}
c_{i,n-i}=\sum\limits_{k=0}^{i}\sum\limits_{l=0}^{n-i}\left( -1\right) ^{i-k}%
\binom{l+i-k}{l}2^{l}r^{l+2i-2k}a_{k,n-i-l}  \label{eq 4}
\end{equation}%
for every $i\in \left\{ 1,\ldots ,n\right\} $. Motivated by expressions (\ref%
{eqpaiva}) and (\ref{eq 4}) define 
\begin{equation*}
\Delta _{i}^{n}\left( r\right) =c_{i,n-i}-\sum\limits_{k=0}^{i}\left(
-1\right) ^{k}\tbinom{n-i+k}{n-i}2^{n-i}r^{n+k}\Gamma
_{i-k}^{n}=\sum\limits_{k=0}^{i}\Delta _{i,k}^{n}\left( r\right)
\end{equation*}%
where%
\begin{eqnarray*}
\Delta _{i,k}^{n}\left( r\right) &=&\sum\limits_{l=0}^{n-i}\left( -1\right)
^{i-k}\binom{l+i-k}{l}2^{l}r^{l+2i-2k}a_{k,n-i-l} \\
&&\mathbb{\qquad }-\sum\limits_{m=0}^{i-k}\sum\limits_{l=0}^{n-i+m}\left(
-1\right) ^{m}\frac{2^{k+m+n}r^{m+n}}{2^{l+2i}r^{l+k}}\tbinom{n-i+m}{n-i}%
\tbinom{i-m-k+l}{l}a_{k,i-m-k+l}\;,
\end{eqnarray*}%
for every $i\in \left\{ 0,\ldots ,n\right\} $ and $k\in \left\{ 0,\ldots
,i\right\} $. Therefore, our goal resumes to verify that%
\begin{equation*}
\Delta _{i,k}^{n}\left( r\right) =0,
\end{equation*}%
for every $i\in \left\{ 0,\ldots ,n\right\} $ and $k\in \left\{ 0,\ldots
,i\right\} $.

\bigskip

Fix $\left( i,k\right) \in \left\{ 0,\ldots ,n\right\} \times \left\{
0,\ldots ,i\right\} $ and define the variables%
\begin{equation*}
p=i-k\in \mathbb{N}\;\text{,}\qquad q=n-i\in \mathbb{N\qquad }\text{and}%
\qquad A_{j}=a_{k,j}\in \mathbb{R}\text{.}
\end{equation*}%
Using these new variables and the identities (\ref{inv_cont_0}) and (\ref%
{insane}) we obtain that%
\begin{equation*}
\Delta _{i,k}^{n}\left( r\right) =\Lambda _{p,q}\left( r\right)
\end{equation*}%
where%
\begin{eqnarray}
\Lambda _{p,q}\left( r\right) &=&\sum\limits_{l=0}^{q}\left( -1\right) ^{p}%
\tbinom{q-l+p}{q-l}\frac{2^{q}r^{q+2p}}{2^{l}r^{l}}A_{l}  \label{Jpq} \\
&&\mathbb{\qquad }-\sum\limits_{l=0}^{p}\left( -1\right) ^{p}\frac{%
2^{q}r^{2p+q}}{2^{l}r^{l}}\left( \sum\limits_{m=0}^{l}\left( -1\right) ^{m}%
\tbinom{q+p-m}{q}\tbinom{l}{m}\right) A_{l}  \notag \\
&&\mathbb{\qquad }-\sum\limits_{l=p+1}^{p+q}\left( -1\right) ^{p}\frac{%
2^{q}r^{2p+q}}{2^{l}r^{l}}\left( \sum\limits_{m=0}^{p}\left( -1\right) ^{m}%
\tbinom{q+p-m}{q}\tbinom{l}{m}\right) A_{l}\;.  \notag
\end{eqnarray}%
Since%
\begin{equation*}
\binom{q+p-m}{q}=0
\end{equation*}%
for every $m>p$, we can rewrite expression (\ref{Jpq}) as%
\begin{eqnarray}
\Lambda _{p,q}\left( r\right) &=&\sum\limits_{l=0}^{q}\left( -1\right) ^{p}%
\frac{2^{q}r^{2p+q}}{2^{l}r^{l}}\left( \tbinom{q-l+p}{q-l}%
-\sum\limits_{m=0}^{l}\left( -1\right) ^{m}\tbinom{q+p-m}{q}\tbinom{l}{m}%
\right) A_{l}  \label{Jpq2} \\
&&-\sum\limits_{l=q+1}^{p+q}\left( -1\right) ^{p}\frac{2^{q}r^{2p+q}}{%
2^{l}r^{l}}\left( \sum\limits_{m=0}^{l}\left( -1\right) ^{m}\tbinom{q+p-m}{q}%
\tbinom{l}{m}\right) A_{l}\;.  \notag
\end{eqnarray}%
Follows from Lemma (\ref{Prop 10.20}) that%
\begin{equation}
\sum_{m=0}^{l}\left( -1\right) ^{m}\dbinom{q+p-m}{q}\dbinom{l}{m}=\dbinom{%
q+p-l}{q-l}  \label{Bozo2}
\end{equation}%
for every $l\in \mathbb{N}$. Combining expressions (\ref{Jpq2}) and (\ref%
{Bozo2}) with the fact that%
\begin{equation*}
\dbinom{q+p-l}{q-l}=0
\end{equation*}%
for every $l>q$, we obtain that $\Lambda _{p,q}\left( r\right) =0$.

\bigskip

\noindent \textbf{Claim 2:} $c_{i,j}=0$, for every $i,j\in \mathbb{N}$ such
that $i+j>n$.

\bigskip

Given $i,j\in \mathbb{N}$ such that $i+j>n$, define $l=i+j-n\in \mathbb{N-}%
\left\{ 0\right\} $. It follows from property (2) and the definition of the
family $\mathcal{A}$ that%
\begin{equation}
c_{i,j}=c_{i,n+l-i}=-r^{2}c_{i-1,n+l-i}+2rc_{i,n+l-i-1}  \label{cij meio}
\end{equation}

We will proceed by induction on $l$. Suppose $l=1$ and observe that (\ref%
{cij meio}) is rewritten as 
\begin{equation*}
c_{i,j}=-r^{2}c_{i-1,n+1-i}+2rc_{i,n-i}\text{.}
\end{equation*}%
First of all, $c_{i,n-i}=0$ by Claim (1). If $i\geq 1$, Claim (1) also
allows us to conclude that $c_{i-1,n+1-i}=0$. On th other hand, we can
conclude the same thing, in case $i=0$, using definition \ref{def cij}.
Therefore, we must have $c_{i,j}=0$.

\bigskip

Assume that $c_{i,n+l-i}=0$ for some $l\geq 1$ and for every $i,j\in \mathbb{%
N}$ such that $l=i+j-n$. Follows from (\ref{cij meio}) that 
\begin{equation*}
c_{i,j}=-r^{2}c_{i-1,n+l+1-i}+2rc_{i,n+l-i}
\end{equation*}%
for every $i,j\in \mathbb{N}$ such that $l+1=i+j-n$. Finally, observe that
induction hypothesis (combined with definition (\ref{def cij}) in the cases $%
i=0$ and $j=0$) allows us to conclude again that $c_{i,j}=0$.

\hfill
\end{proof}

\bigskip

Finally, we finish the section with another technical result that will be
important in the next section to deal with the Lorentzian case.

\begin{theorem}
\label{Q_Qepsilon}Given $\varepsilon \in \left\{ -1,1\right\} $ and a non
constant polynomial%
\begin{equation*}
Q\left( x,y\right)
=\sum\limits_{i=0}^{n}\sum\limits_{j=0}^{n-i}a_{i,j}x^{i}y^{j}\in \mathbb{R}%
\left[ x,y\right] \text{,}
\end{equation*}%
define the polynomial%
\begin{equation*}
Q_{\varepsilon }\left( x,y\right)
=\sum\limits_{i=0}^{n}\sum\limits_{j=0}^{n-i}\varepsilon
^{i+j}a_{i,j}x^{i}y^{j}\in \mathbb{R}\left[ x,y\right] .
\end{equation*}%
For every nonzero constant $r\in \mathbb{R}$, we have that $Q_{\varepsilon
}\left( x,y\right) $ belongs to the ideal generated by the polynomial $%
xr^{2}-2ry+1$ if and only if $Q\left( x,y\right) $ belongs to the ideal
generated by the polynomial $xr^{2}-2ry+\varepsilon $.
\end{theorem}

\begin{proof}
Suppose that $Q_{\varepsilon }\left( x,y\right) $ belongs to the ideal
generated by $xr^{2}-2ry+1$. Then there is a polynomial%
\begin{equation*}
R\left( x,y\right)
=\sum\limits_{i=0}^{n-1}\sum\limits_{j=0}^{n-1-i}b_{i,j}x^{i}y^{j}
\end{equation*}%
such that%
\begin{equation*}
Q_{\varepsilon }\left( x,y\right) =\left( xr^{2}-2ry+1\right) R\left(
x,y\right) \text{.}
\end{equation*}%
Defining $b_{i,j}=0$, whenever $i<0$, $j<0$ or $i+j>n-1$, we have%
\begin{eqnarray*}
\sum\limits_{i=0}^{n}\sum\limits_{j=0}^{n-i}\varepsilon
^{i+j}a_{i,j}x^{i}y^{j} &=&\left( xr^{2}-2ry+1\right)
\sum\limits_{i=0}^{n-1}\sum\limits_{j=0}^{n-1-i}b_{i,j}x^{i}y^{j} \\
&=&\sum\limits_{i=0}^{n-1}\sum\limits_{j=0}^{n-1-i}r^{2}b_{i,j}x^{i+1}y^{j}-%
\sum\limits_{i=0}^{n-1}\sum\limits_{j=0}^{n-1-i}2rb_{i,j}x^{i}y^{j+1} \\
&&\mathbb{\qquad }+\sum\limits_{i=0}^{n-1}\sum%
\limits_{j=0}^{n-1-i}b_{i,j}x^{i}y^{j} \\
&=&\sum\limits_{i=0}^{n}\sum\limits_{j=0}^{n-i}\left(
r^{2}b_{i-1,j}-2rb_{i,j-1}+b_{i,j}\right) x^{i}y^{j}
\end{eqnarray*}%
and consequentely that%
\begin{equation}
a_{i,j}=\frac{r^{2}b_{i-1,j}-2rb_{i,j-1}+b_{i,j}}{\varepsilon ^{i+j}}\text{,}
\label{a_ij}
\end{equation}%
for every $i\in \left\{ 0,\ldots ,n\right\} $ and $j\in \left\{ 0,\ldots
,n-i\right\} $.

\bigskip

To prove that $Q\left( x,y\right) $ belongs to the ideal generated by the
polynomial $xr^{2}-2ry+\varepsilon $ is sufficient to prove that the
polynomial%
\begin{equation*}
T\left( x,y\right)
=\sum\limits_{i=0}^{n}\sum\limits_{j=0}^{n-i}a_{i,j}x^{i}y^{j}-\left(
r^{2}x-2ry+\varepsilon \right)
\sum\limits_{i=0}^{n-1}\sum\limits_{j=0}^{n-1-i}\varepsilon
^{-i-j-1}b_{i,j}x^{i}y^{j}\in \mathbb{R}\left[ x,y\right]
\end{equation*}%
is identicaly null. Using equality (\ref{a_ij}) we obtain%
\begin{eqnarray*}
T\left( x,y\right) &=&\sum\limits_{i=0}^{n}\sum\limits_{j=0}^{n-i}\frac{%
r^{2}b_{i-1,j}}{\varepsilon ^{i+j}}x^{i}y^{j}-\sum\limits_{i=0}^{n}\sum%
\limits_{j=0}^{n-i}\frac{2rb_{i,j-1}}{\varepsilon ^{i+j}}x^{i}y^{j} \\
&&\mathbb{\qquad }+\sum\limits_{i=0}^{n}\sum\limits_{j=0}^{n-i}\frac{b_{i,j}%
}{\varepsilon ^{i+j}}x^{i}y^{j}-\sum\limits_{i=0}^{n-1}\sum%
\limits_{j=0}^{n-1-i}r^{2}\varepsilon ^{-i-j-1}b_{i,j}x^{i+1}y^{j} \\
&&\mathbb{\qquad }+\sum\limits_{i=0}^{n-1}\sum\limits_{j=0}^{n-1-i}2r%
\varepsilon
^{-i-j-1}b_{i,j}x^{i}y^{j+1}-\sum\limits_{i=0}^{n-1}\sum%
\limits_{j=0}^{n-1-i}\varepsilon ^{-i-j}b_{i,j}x^{i}y^{j}
\end{eqnarray*}%
The nullity of the polynomial $T\left( x,y\right) $ is a consequence of the
following equalities whose veracity derives from the definition of the
coeficients $b_{i,j}$:%
\begin{equation*}
\sum\limits_{i=0}^{n}\sum\limits_{j=0}^{n-i}\frac{r^{2}b_{i-1,j}}{%
\varepsilon ^{i+j}}x^{i}y^{j}=\sum\limits_{i=0}^{n-1}\sum%
\limits_{j=0}^{n-1-i}r^{2}\varepsilon ^{-i-j-1}b_{i,j}x^{i+1}y^{j}
\end{equation*}%
\begin{equation*}
\sum\limits_{i=0}^{n}\sum\limits_{j=0}^{n-i}\frac{2rb_{i,j-1}}{\varepsilon
^{i+j}}x^{i}y^{j}=\sum\limits_{i=0}^{n-1}\sum\limits_{j=0}^{n-1-i}2r%
\varepsilon ^{-i-j-1}b_{i,j}x^{i}y^{j+1}
\end{equation*}%
\begin{equation*}
\sum\limits_{i=0}^{n}\sum\limits_{j=0}^{n-i}\frac{b_{i,j}}{\varepsilon ^{i+j}%
}x^{i}y^{j}=\sum\limits_{i=0}^{n-1}\sum\limits_{j=0}^{n-1-i}\frac{b_{i,j}}{%
\varepsilon ^{i+j}}x^{i}y^{j}.
\end{equation*}

\bigskip

Suppose now that $Q\left( x,y\right) $ belongs to the ideal generated by $%
xr^{2}-2ry+\varepsilon $. Then there is a polynomial

\begin{equation*}
U\left( x,y\right)
=\sum\limits_{i=0}^{n-1}\sum\limits_{j=0}^{n-1-i}b_{i,j}x^{i}y^{j}
\end{equation*}%
such that%
\begin{equation*}
Q\left( x,y\right) =\left( xr^{2}-2ry+\varepsilon \right) U\left( x,y\right) 
\text{.}
\end{equation*}%
Following the same steps we conclude that%
\begin{equation*}
Q_{\varepsilon }\left( x,y\right) =\left( r^{2}x-2ry+1\right)
\sum\limits_{i=0}^{n-1}\sum\limits_{j=0}^{n-1-i}\varepsilon
^{i+j+1}c_{i,j}x^{i}y^{j}\in \left\langle r^{2}x-2ry+1\right\rangle .
\end{equation*}

\hfill
\end{proof}

\begin{corollary}
\label{Teo_Factorization_epsilon}Consider $\varepsilon \in \left\{
-1,1\right\} $, a non constant polynomial $Q\left( x,y\right) $ and a
nonzero constant $r\in \mathbb{R}$. Then $Q\left( x,y\right) $ belongs to
the ideal generated by $xr^{2}-2yr+\varepsilon $ if and only if the
polynomial $Q\left( \frac{\varepsilon x}{r},\frac{\varepsilon \left(
xr+1\right) }{2r}\right) \in \mathbb{R}\left[ x\right] $ is identically null.
\end{corollary}

\begin{proof}
Write $Q\left( x,y\right)
=\sum\limits_{i=0}^{n}\sum\limits_{j=0}^{n-i}a_{i,j}x^{i}y^{j}\in \mathbb{R}%
\left[ x,y\right] $ and define the polynomial%
\begin{equation*}
Q_{\varepsilon }\left( x,y\right)
=\sum\limits_{i=0}^{n}\sum\limits_{j=0}^{n-i}\varepsilon
^{i+j}a_{i,j}x^{i}y^{j}\in \mathbb{R}\left[ x,y\right] .
\end{equation*}%
Since%
\begin{eqnarray*}
Q_{\varepsilon }\left( \frac{x}{r},\frac{xr+1}{2r}\right)
&=&\sum\limits_{i=0}^{n}\sum\limits_{j=0}^{n-i}\varepsilon
^{i+j}a_{i,j}\left( \frac{x}{r}\right) ^{i}\left( \frac{xr+1}{2r}\right) ^{j}
\\
&=&\sum\limits_{i=0}^{n}\sum\limits_{j=0}^{n-i}a_{i,j}\left( \frac{%
\varepsilon x}{r}\right) ^{i}\left( \frac{\varepsilon \left( xr+1\right) }{2r%
}\right) ^{j} \\
&=&Q\left( \frac{\varepsilon x}{r},\frac{\varepsilon \left( xr+1\right) }{2r}%
\right) ,
\end{eqnarray*}%
the result follows from Theorems \ref{Teo_Factorization} and \ref{Q_Qepsilon}%
.
\end{proof}

\section{Main result and applications for Tubular surfaces\label{main result
and applications}}

In this section we present our main theorems that fully classify Polynomial
Weingarten tubular surfaces in Euclidean, Lorentzian and Hyperbolic $3$%
-space. Our presentation will focus on the Lorentzian case. It will be clear
from the proofs that the Euclidean and Hiperbolic theorems are analogous and
less technical versions of the Lorentzian ones. We start with some algebraic
definitions:

\bigskip

Given $\varepsilon \in \left\{ -1,1\right\} $, the \textbf{Lorentzian }$%
\varepsilon $\textbf{-radius} of a polynomial $Q\left( x,y\right) \in 
\mathbb{R}\left[ x,y\right] $ is defined as the set%
\begin{equation*}
Rad_{L,\varepsilon }\left( Q\right) =\left\{ r\in \left( 0,+\infty \right) 
\text{ };\text{ }Q\left( 0,\tfrac{\varepsilon }{2r}\right) =0\right\} \text{.%
}
\end{equation*}%
The \textbf{Lorentzian }$\varepsilon $\textbf{-star radius} of $Q\left(
x,y\right) $ is defined as the set 
\begin{equation*}
Rad_{L,\varepsilon }^{\ast }\left( Q\right) =\left\{ r\in Rad_{L,\varepsilon
}\left( Q\right) \;;\;Q\left( x,y\right) \in \left\langle
xr^{2}-2ry+\varepsilon \right\rangle \right\} \text{,}
\end{equation*}%
where $\left\langle \varphi \left( x,y\right) \right\rangle $ denotes the
ideal in $\mathbb{R}\left[ x,y\right] $ generated by the polynomial $\varphi
\left( x,y\right) $. The \textbf{Lorentzian radius} and the \textbf{%
Lorentzian star radius }of $Q\left( x,y\right) $ are respectively defined by%
\begin{equation*}
Rad_{L}\left( Q\right) =\bigcup\limits_{\varepsilon \in \left\{ -1,1\right\}
}Rad_{L,\varepsilon }\left( Q\right) \mathbb{\qquad }\text{and}\qquad
Rad_{L}^{\ast }\left( Q\right) =\bigcup\limits_{\varepsilon \in \left\{
-1,1\right\} }Rad_{L,\varepsilon }^{\ast }\left( Q\right) \text{.}
\end{equation*}

\begin{example}
\label{Example3}Consider a polynomial $Q\left( x,y\right) \in \mathbb{R}%
\left[ x,y\right] $and real number $r>0$. By Proposition \ref{Prop
Polynomial 1} we have%
\begin{equation*}
Q\left( \tfrac{x}{r},\tfrac{xr+1}{2r}\right) =\sum\limits_{k=0}^{n}\Gamma
_{k}\left( r\right) x^{k}\text{,}
\end{equation*}%
where the coefficients $\Gamma _{k}\left( r\right) $ are defined in (\ref%
{def Tau}). Then, we have $r\in Rad_{L,1}\left( Q\right) $ if and only if $%
\tfrac{1}{2r}$ is a positive root of the polynomial $\Gamma _{0}$. In
particular, it follows from Example (\ref{Exemplo2}) that the polynomial $%
Q\left( x,y\right) $ presented in the Example \ref{Ex Q} verify%
\begin{equation*}
Rad_{L,1}\left( Q\right) =\left\{ \frac{1}{8},2\right\} .
\end{equation*}
\end{example}

\bigskip

Now we present the Lorentzian version of our first main theorem:

\begin{theorem}[Lorentzian $\mathcal{S}\left( Q\right) $ problem]
\label{Teo S(Q) Lorentzian}Given a polynomial $Q\left( x,y\right) \in 
\mathbb{R}\left[ x,y\right] $, denote by $\mathcal{S}_{L}\left( Q\right) $
the set of all non degenerated regular tubular surfaces in Lorentzian $3$%
-space verifying $Q\left( K,H\right) \equiv 0$. Then, the elements of $%
\mathcal{S}_{L}\left( Q\right) $ are:

\begin{enumerate}
\item[i.] The right cylinders with signal $\varepsilon \in \left\{
-1,1\right\} $ and radius in $Rad_{L,\varepsilon }\left( Q\right) $;

\item[ii.] The regular tubular surfaces with signal $\varepsilon \in \left\{
-1,1\right\} $ and radius in $Rad_{L,\varepsilon }^{\ast }\left( Q\right) $.
\end{enumerate}
\end{theorem}

\begin{proof}
Having zero Gaussian curvature, it follows directly from the definition of $%
Rad_{L,\varepsilon }\left( Q\right) $ that right cylinders with signal $%
\varepsilon \in \left\{ -1,1\right\} $ and radius in $Rad_{L,\varepsilon
}\left( Q\right) $ belong to $\mathcal{S}_{L}\left( Q\right) $.

Take a regular tubular surface $S$ with signal $\varepsilon \in \left\{
-1,1\right\} $ and with radius in $Rad_{L,\varepsilon }^{\ast }\left(
Q\right) $ (assuming it is not empty). Follows from the definition of $%
Rad_{L,\varepsilon }^{\ast }\left( Q\right) $ that there is a polynomial $%
R\left( x,y\right) \in \mathbb{R}\left[ x,y\right] $ such that%
\begin{equation*}
Q\left( x,y\right) =\left( xr^{2}-2ry+\varepsilon \right) R\left( x,y\right) 
\text{.}
\end{equation*}%
Using the expressions for the curvatures of $S$ presented in (\ref%
{Tubular_curvatures}) we obtain%
\begin{eqnarray*}
Q\left( K,H\right) &=&\left( Kr^{2}-2rH+\varepsilon \right) R\left(
K,H\right) \\
&=&\left( \varepsilon \dfrac{\varepsilon _{B}\kappa \mu }{r\left(
1+\varepsilon _{B}r\kappa \mu \right) }r^{2}-2r\varepsilon \dfrac{%
2\varepsilon _{B}r\kappa \mu +1}{2r\left( 1+\varepsilon _{B}r\kappa \mu
\right) }+\varepsilon \right) R\left( K,H\right) \\
&=&\left( \dfrac{\varepsilon \varepsilon _{B}\kappa \mu r}{1+\varepsilon
_{B}r\kappa \mu }-\dfrac{2\varepsilon \varepsilon _{B}r\kappa \mu
+\varepsilon }{1+\varepsilon _{B}r\kappa \mu }+\frac{\varepsilon
+\varepsilon \varepsilon _{B}r\kappa \mu }{1+\varepsilon _{B}r\kappa \mu }%
\right) R\left( K,H\right) \\
&=&0\text{.}
\end{eqnarray*}%
The equality above shows that $S\in \mathcal{S}_{L}\left( Q\right) $.

\bigskip

Now take in $\mathcal{S}_{L}\left( Q\right) $ (assuming it is not empty) a
non degenerated regular tubular surface $S$ with signal $\varepsilon \in
\left\{ -1,1\right\} $ and radius $r>0$. Follows from expressions (\ref%
{Tubular_curvatures}) that the Gaussian and the mean curvatures of $S$ are%
\begin{equation*}
K=\varepsilon k_{1}k_{2}=\mathbb{\qquad }\text{and}\qquad H=\varepsilon 
\frac{k_{1}+k_{2}}{2}\text{,}
\end{equation*}%
where%
\begin{equation*}
k_{1}=\dfrac{\varepsilon _{B}\kappa \mu }{1+\varepsilon _{B}r\kappa \mu }%
\mathbb{\qquad }\text{and}\qquad k_{2}=\frac{1}{r}\text{.}
\end{equation*}%
It is important to remark that the principal curvatures of a regular surface
in the Lorentzian $3$-space are not always well defined. In many cases, the
functions $k_{1}$ and $k_{2}$ defined above are the principal curvatures of $%
S$, but, in general, they are just functions defined on $S$.

\bigskip

Write $Q\left( x,y\right)
=\sum\limits_{i=0}^{n}\sum\limits_{j=0}^{n-i}a_{i,j}x^{i}y^{j}$ and define
the polynomial%
\begin{equation*}
Q_{\varepsilon }\left( x,y\right)
=\sum\limits_{i=0}^{n}\sum\limits_{j=0}^{n-i}\varepsilon
^{i+j}a_{i,j}x^{i}y^{j}\in \mathbb{R}\left[ x,y\right] .
\end{equation*}%
Note that%
\begin{eqnarray*}
0 &=&Q\left( K,H\right) =Q\left( \varepsilon k_{1}k_{2},\varepsilon \frac{%
k_{1}+k_{2}}{2}\right)
=\sum\limits_{i=0}^{n}\sum\limits_{j=0}^{n-i}a_{i,j}\left( \varepsilon
k_{1}k_{2}\right) ^{i}\left( \varepsilon \frac{k_{1}+k_{2}}{2}\right) ^{j} \\
&=&\sum\limits_{i=0}^{n}\sum\limits_{j=0}^{n-i}\varepsilon
^{i+j}a_{i,j}\left( \frac{k_{1}}{r}\right) ^{i}\left( \frac{rk_{1}+1}{2r}%
\right) ^{j}=Q_{\varepsilon }\left( \frac{k_{1}}{r},\frac{rk_{1}+1}{2r}%
\right) .
\end{eqnarray*}%
Therefore, $k_{1}\left( s,t\right) $ is a root of $Q_{\varepsilon }\left( 
\frac{x}{r},\frac{xr+1}{2r}\right) \in \mathbb{R}\left[ x\right] $, for
every $\left( s,t\right) \in \left( a,b\right) \times \mathbb{R}$.

\bigskip

If $Q_{\varepsilon }\left( \frac{x}{r},\frac{xr+1}{2r}\right) $ is not the
null polynomial, we have that $k_{1}\left( s,t\right) $ is a constant
function. That is, there is $\lambda \in \mathbb{R}$ such that%
\begin{equation*}
\dfrac{\kappa \left( s\right) \mu \left( t\right) \varepsilon _{B}}{%
1+r\kappa \left( s\right) \mu \left( t\right) \varepsilon _{B}}=\lambda 
\text{,}
\end{equation*}%
for every $\left( s,t\right) \in \left( a,b\right) \times \mathbb{R}$. If $%
\lambda \neq 0$, we can write%
\begin{equation*}
\mu \left( t\right) =\frac{\lambda }{\varepsilon _{B}\kappa \left( s\right)
\left( 1-r\lambda \right) }\text{,}
\end{equation*}%
for every $\left( s,t\right) \in \left( a,b\right) \times \mathbb{R}$. The
above equality implies th $\mu \left( t\right) $ is constant, which is an
absurd because $\mu \left( t\right) $ represents a trigonometric function.
Then, we must have $\lambda =0$ and this fact implies that $r\in
Rad_{L,\varepsilon }\left( Q\right) $ because%
\begin{equation*}
0=Q\left( K,H\right) =Q\left( \varepsilon k_{1}k_{2},\varepsilon \frac{%
k_{1}+k_{2}}{2}\right) =Q\left( 0,\frac{\varepsilon }{2r}\right) \text{.}
\end{equation*}%
More over, since $\mu $ is not null, we can also conclude that $\kappa
\equiv 0$. Therefore $S$ is a right cylinder.

\bigskip

Now let us deal with the case where $Q_{\varepsilon }\left( \frac{x}{r},%
\frac{xr+1}{2r}\right) $ is the null polynomial. It follows from Theorem (%
\ref{Teo_Factorization}) that $Q_{\varepsilon }\left( x,y\right) $ belongs
to the ideal generated by the polynomial $xr^{2}-2yr+1\in \mathbb{R}\left[
x,y\right] $ and from Theorem (\ref{Q_Qepsilon}) that $Q\left( x,y\right) $
belongs to the ideal generated by the polynomial $xr^{2}-2ry+\varepsilon \in 
\mathbb{R}\left[ x,y\right] $. Therefore, there is a polynomial $R\left(
x,y\right) \in \mathbb{R}\left[ x,y\right] $ such that%
\begin{equation*}
Q\left( x,y\right) =\left( xr^{2}-2ry+\varepsilon \right) R\left( x,y\right) 
\text{.}
\end{equation*}%
To prove that $r\in Rad_{L,\varepsilon }^{\ast }\left( Q\right) $ it is
sufficient (by definition) to show that $r\in Rad_{L,\varepsilon }\left(
Q\right) $, but this is an immediate consequence of the equality above. In
fact%
\begin{equation*}
Q\left( 0,\frac{\varepsilon }{2r}\right) =\left( 0r^{2}-2r\frac{\varepsilon 
}{2r}+\varepsilon \right) R\left( 0,\frac{\varepsilon }{2r}\right) =0.
\end{equation*}

\hfill
\end{proof}

\bigskip

The following corollary is an immediate consequence of the above theorem:

\begin{corollary}[Lorentzian]
Consider a polynomial $Q\left( x,y\right) \in \mathbb{R}\left[ x,y\right] $.
There is a non degenerated regular tubular surfaces in Lorentzian $3$-space
verifying $Q\left( K,H\right) \equiv 0$ if and only if $Rad_{L}\left(
Q\right) \neq \emptyset $. Moreover, the radii of all such surfaces belongs
to $Rad_{L}\left( Q\right) $.
\end{corollary}

\bigskip

Since all regular surfaces in Euclidean and Hyperbolic spaces has signal $+1$%
, the definition of radius and star radius in these spaces become simpler.
As seen in the introduction (see Theorem (\ref{TeoTubIntro}), the Euclidean
definitions for the radius and the star radius of a polynomial $Q\left(
x,y\right) \in \mathbb{R}\left[ x,y\right] $ are:

\begin{equation*}
Rad_{E}\left( Q\right) =\left\{ r\in \left( 0,+\infty \right) \text{ };\text{
}Q\left( 0,\tfrac{1}{2r}\right) =0\right\}
\end{equation*}%
\begin{equation*}
Rad_{E}^{\ast }\left( Q\right) =\left\{ r\in Rad_{E}\left( Q\right)
\,;\,Q\left( x,y\right) \in \left\langle xr^{2}-2ry+1\right\rangle \right\} 
\text{.}
\end{equation*}%
The hyperbolic definitions are slightly more technical, but they are
motivated by the expressions for the curvatures of regular tubular surfaces
in hyperbolic $3$-space presented at the end of the Section \ref{Tubular
Surface}.

\bigskip

Given a polynomial $Q\left( x,y\right) \in \mathbb{R}\left[ x,y\right] $,
the \textbf{Hyperbolic} \textbf{radius} and the \textbf{Hyperbolic star
radius} of $Q\left( x,y\right) $ are defined respectively by%
\begin{equation*}
Rad_{H}\left( Q\right) =\left\{ r\in \left( 0,+\infty \right) \text{ };\text{
}Q\left( 0,\tfrac{1}{2\sinh r}\right) =0\right\} 
\end{equation*}%
\begin{equation*}
Rad_{H}^{\ast }\left( Q\right) =\left\{ r\in Rad_{H}\left( Q\right)
\,;\,Q\left( x,y\right) \in \left\langle x\sinh ^{2}r-2y\left( \sinh
r\right) +1\right\rangle \right\} \text{.}
\end{equation*}

The Euclidean results were stated in the introduction. We present below only
the Hyperbolic versions:

\begin{theorem}[Hyperbolic $\mathcal{S}\left( Q\right) $ problem]
Given a polynomial $Q\left( x,y\right) \in \mathbb{R}\left[ x,y\right] $,
denote by $\mathcal{S}_{H}\left( Q\right) $ the set of all regular tubular
surfaces in Hyperbolic $3$-space verifying $Q\left( K,H\right) \equiv 0$.
Then, the elements of $\mathcal{S}_{H}\left( Q\right) $ are:

\begin{enumerate}
\item[i.] The right cylinders with radius in $Rad_{H}\left( Q\right) $;

\item[ii.] The regular tubular surfaces with radius in $Rad_{H}^{\ast
}\left( Q\right) $.
\end{enumerate}
\end{theorem}

\begin{corollary}[Hyperbolic]
Consider a polynomial $Q\left( x,y\right) \in \mathbb{R}\left[ x,y\right] $.
There is a regular tubular surfaces in Hyperbolic $3$-space verifying $%
Q\left( K,H\right) \equiv 0$ if and only if $Rad_{H}\left( Q\right) \neq
\emptyset $. Moreover, the radii of all such surfaces belongs to $%
Rad_{H}\left( Q\right) $.
\end{corollary}

\bigskip

The example below elucidates the usefulness of the above theorems (together
with Theorem \ref{Teo_Factorization}) with a numerical case:

\begin{example}
\label{Example SQ}Consider the polynomial $Q\left( x,y\right)
=14y-25x+100xy-40y^{2}-1\in \mathbb{R}\left[ x,y\right] $. Since 
\begin{equation*}
Q\left( 0,\tfrac{1}{2r}\right) =-\frac{\left( r-2\right) \left( r-5\right) }{%
r^{2}},
\end{equation*}%
we have $Rad_{E}\left( Q\right) =\left\{ 2,5\right\} $. Simple calculations
allow us to conclude that $Q\left( \frac{x}{5},\frac{5x+1}{10}\right) $
vanishes identically and $Q\left( \frac{x}{2},\frac{2x+1}{4}\right)
=\allowbreak 3x\left( 5x-1\right) $. So, the elements of $\mathcal{S}%
_{E}\left( Q\right) $ are arbitrary regular tubular surfaces of radius $5$
and right cylinders of radius $2$.
\end{example}

\bigskip 

Observe that the above theorems allow us to classify tubular Polynomial
Weingarten surfaces through the study of polynomials (exclusively). To
portray its funcionality, let us show how some well-known results can be
obtained from it:

\begin{example}
Given $c\in \mathbb{R}$, consider the polynomial $Q\left( x,y\right) =y-c\in 
\mathbb{R}\left[ x,y\right] $ and note that $Rad_{E}\left( Q\right) \neq
\emptyset $ if and only if $c>0$. More precisely,%
\begin{equation*}
Rad_{E}\left( Q\right) =\left\{ 
\begin{array}{ll}
\left\{ \frac{1}{2c}\right\} & \text{if }c>0 \\ 
\emptyset & \text{if }c\leq 0%
\end{array}%
\right. .
\end{equation*}%
Then, there are no regular tubular surfaces with negative constant mean
curvature. Moreover, since 
\begin{equation*}
Q\left( \frac{x}{\frac{1}{2c}},\frac{x\frac{1}{2c}+1}{2\frac{1}{2c}}\right)
=Q\left( 2cx,\frac{x+2c}{2}\right) =\frac{x+2c}{2}-c=\frac{x}{2}
\end{equation*}%
is not null, the unique regular tubular surfaces with constant mean
curvature $c>0$ are the right cylinders of radius $\frac{1}{2c}$.
\end{example}

\begin{example}
Given $c\in \mathbb{R}$, consider the polynomial $Q\left( x,y\right) =x-c\in 
\mathbb{R}\left[ x,y\right] $ and note that $Rad_{E}\left( Q\right) \neq
\emptyset $ if and only if $c=0$. More precisely,%
\begin{equation*}
Rad_{E}\left( Q\right) =\left\{ 
\begin{array}{ll}
\left( 0,+\infty \right) & \text{if }c=0 \\ 
\emptyset & \text{if }c\neq 0%
\end{array}%
\right. .
\end{equation*}%
When $c=0$ we have that 
\begin{equation*}
Q\left( \tfrac{x}{r},\tfrac{xr+1}{2r}\right) =\tfrac{x}{r}
\end{equation*}%
is not null, for every $r\in \left( 0,+\infty \right) $. Then The cylinders
(of any radii) are the only tubular surfaces with constant Gaussian
curvature.
\end{example}

\bigskip

As mentioned in the introduction, a special and very studied polynomial
relation among the Gaussian and the mean curvatures, is the linear one. The
classification of regular tubular surfaces in Euclidean $3$-space whose
Gaussian and mean curvatures verify a linear relation was stated in the
introduction. Let us present below the Lorentzian and hyperbolic cases. We
need the following notation:%
\begin{equation*}
\func{sgn}\left( \lambda \right) =\left\{ 
\begin{array}{ll}
+1 & \text{if }\lambda >0 \\ 
0 & \text{if }\lambda =0 \\ 
-1 & \text{if }\lambda <0%
\end{array}%
\right.
\end{equation*}%
for every $\lambda \in \mathbb{R}$.

\begin{corollary}[Lorentzian]
Given real numbers $a,b,c$ such that $\left( a,b\right) \neq \left(
0,0\right) $, consider the polynomial 
\begin{equation*}
Q\left( x,y\right) =ax+by-c\in \mathbb{R}\left[ x,y\right]
\end{equation*}%
Then, $\mathcal{S}_{L}\left( Q\right) \neq \emptyset $ if and only if $b,c=0$
or $b,c\neq 0$. More precisely, for $\varepsilon =\func{sgn}\left( bc\right) 
$ and $\Delta =\varepsilon b^{2}+4ac$, we have:

\begin{enumerate}
\item[i.] If $b,c=0$, then $\mathcal{S}_{L}\left( Q\right) $ is the set of
all right cylinders of any signal and any radius;

\item[ii.] If $b,c\neq 0$ and $\Delta =0$, then $\mathcal{S}_{L}\left(
Q\right) $ is the set of all regular tubular surfaces with signal $%
\varepsilon $ and radius $\frac{b\varepsilon }{2c}$;

\item[iii.] If $b,c\neq 0$ and $\Delta \neq 0$, then $\mathcal{S}_{L}\left(
Q\right) $ is the set of all right cylinders with signal $\varepsilon $ and
radius $\frac{b\varepsilon }{2c}$.
\end{enumerate}
\end{corollary}

\begin{proof}
Since $Q\left( 0,\tfrac{\varepsilon }{2r}\right) =b\tfrac{\varepsilon }{2r}%
-c $, an easy calculation give us%
\begin{equation*}
Rad_{L,\varepsilon }\left( Q\right) =\left\{ 
\begin{array}{ll}
\left( 0,+\infty \right) & \text{if }b=c=0 \\ 
\tfrac{b\varepsilon }{2c} & \text{if }\varepsilon bc>0 \\ 
\emptyset & \text{otherwise}%
\end{array}%
\right. .
\end{equation*}%
First of all, observe that%
\begin{equation*}
Q\left( \frac{\varepsilon x}{r},\frac{\varepsilon \left( xr+1\right) }{2r}%
\right) =a\frac{\varepsilon x}{r}+b\frac{\varepsilon \left( xr+1\right) }{2r}%
-c=\varepsilon \frac{2a+br}{2r}x+\frac{b\varepsilon -2cr}{2r}.
\end{equation*}%
When $b=c=0$, we have $a\neq 0$ by hypothesis and consequently%
\begin{equation*}
Q\left( \frac{\varepsilon x}{r},\frac{\varepsilon \left( xr+1\right) }{2r}%
\right) =\frac{a\varepsilon }{r}x\not\equiv 0,
\end{equation*}%
for every $r\in \left( 0,+\infty \right) $. Suppose now the case $%
\varepsilon bc>0$ (which implies that $b,c\neq 0$). For $r=\tfrac{%
b\varepsilon }{2c}$ we have%
\begin{equation*}
Q\left( \frac{\varepsilon x}{r},\frac{\varepsilon \left( xr+1\right) }{2r}%
\right) =\varepsilon \frac{2a+b\left( \tfrac{b\varepsilon }{2c}\right) }{%
2\left( \tfrac{b\varepsilon }{2c}\right) }x+\frac{b\varepsilon -2c\left( 
\tfrac{b\varepsilon }{2c}\right) }{2\left( \tfrac{b\varepsilon }{2c}\right) }%
=\frac{\Delta }{2b}x\text{.}
\end{equation*}%
Then, for $r=\tfrac{b\varepsilon }{2c}$, the polynomial $Q\left( \frac{%
\varepsilon x}{r},\frac{\varepsilon \left( xr+1\right) }{2r}\right) $ is
identically null if and only if the discriminant $\Delta $ is zero. To
conclude what was stated in this corollary, we must only unite all the above
facts with Theorem \ref{Teo S(Q) Lorentzian} and Corollary \ref%
{Teo_Factorization_epsilon}.\hfill
\end{proof}

\begin{corollary}[Hyperbolic]
Given real numbers $a,b,c$ such that $\left( a,b\right) \neq \left(
0,0\right) $, consider the polynomial%
\begin{equation*}
Q\left( x,y\right) =ax+by-c\in \mathbb{R}\left[ x,y\right] \text{.}
\end{equation*}%
Then, $\mathcal{S}_{H}\left( Q\right) \neq \emptyset $ if and only if $b,c=0$
or $bc>0$. More precisely, for $\Delta =b^{2}+4ac$, we have:

\begin{enumerate}
\item[i.] If $b,c=0$, then $\mathcal{S}_{H}\left( Q\right) $ is the set of
all right cylinders of any radius;

\item[ii.] If $bc>0$ and $\Delta =0$, then $\mathcal{S}_{H}\left( Q\right) $
is the set of all tubular surfaces of radius $\sinh ^{-1}\left( \frac{b}{2c}%
\right) $;

\item[iii.] If $bc>0$ and $\Delta \neq 0$, then $\mathcal{S}_{H}\left(
Q\right) $ is the set of all right cylinders of radius $\sinh ^{-1}\left( 
\frac{b}{2c}\right) $.
\end{enumerate}
\end{corollary}

\bigskip

The length of the second fundamental form of a surface is defined by%
\begin{equation*}
\left\vert \mathcal{A}\right\vert =\sqrt{4H^{2}-2K}.
\end{equation*}%
Then surfaces with second fundamental form of constant length, say $%
\left\vert \mathcal{A}\right\vert =c>0$, are Weingarten surfaces whose
Gaussian and mean curvatures vanish the polynomial%
\begin{equation*}
Q_{c}\left( x,y\right) =-2x+4y^{2}-c^{2}.
\end{equation*}

Again we present below just the Lorentzian and the Hyperbolic
classifications. The Euclidean on was stated in the introduction.

\begin{corollary}[Lorentzian]
The right cylinders are the unique non degenerated regular tubular surfaces
in th Lorentzian $3$-space with second fundamental form of constant length.
More precisely, the unique non degenerated regular tubular surfaces in the
Lorentzian $3$-space verifying $\left\vert \mathcal{A}\right\vert =c>0$ are
the right cylinders of radius $\frac{1}{c}$.
\end{corollary}

\begin{proof}
Given $c>0$ and $\varepsilon \in \left\{ -1,1\right\} $, we have%
\begin{equation*}
Q_{c}\left( 0,\tfrac{\varepsilon }{2r}\right) =4\left( \tfrac{\varepsilon }{%
2r}\right) ^{2}-c^{2}=\frac{1-c^{2}r^{2}}{r^{2}}
\end{equation*}%
and therefore%
\begin{equation*}
Rad_{L,\varepsilon }\left( Q_{c}\right) =\left\{ \frac{1}{c}\right\} .
\end{equation*}%
Since%
\begin{eqnarray*}
Q_{c}\left( \frac{\varepsilon x}{\frac{1}{c}},\frac{\varepsilon \left( x%
\frac{1}{c}+1\right) }{2\frac{1}{c}}\right) &=&Q_{c}\left( cx\varepsilon ,%
\frac{\varepsilon \left( x+c\right) }{2}\right) \\
&=&-2cx\varepsilon +4\left( \frac{\varepsilon \left( x+c\right) }{2}\right)
^{2}-c^{2} \\
&=&x^{2}+2c\left( 1-\varepsilon \right) x\not\equiv 0
\end{eqnarray*}%
the statement of the corollary is a consequence of , we must only unite all
the above facts with he Theorem \ref{Teo S(Q) Lorentzian} and the Corollary %
\ref{Teo_Factorization_epsilon}.\hfill
\end{proof}

\begin{corollary}[Hyperbolic]
The right cylinders are the unique non degenerated regular tubular surfaces
in th Hyperbolic $3$-space with second fundamental form of constant length.
More precisely, the unique non degenerated regular tubular surfaces in the
Hi[erbolic $3$-space verifying $\left\vert \mathcal{A}\right\vert =c>0$ are
the right cylinders of radius $\sinh ^{-1}\left( \frac{1}{c}\right) $.
\end{corollary}

\bigskip

The remaining question will be answered by the our last main Theorem, once
it acts in the another direction of the previous one. More precisely, now we
consider a (fixed) regular tubular surface $S$ and we present the set $%
\mathcal{Q}\left( S\right) $ of all polynomials that vanishes its
curvatures. Furthermore, an algebraic characterizations of these polynomials
are expressed as a discriminant to analyse tubular Polynomial Weingarten
surfaces.

\begin{theorem}[Lorentzian $\mathcal{Q}\left( S\right) $ problem]
Given a non degenerated regular tubular surface $S$ of radius $r>0$ and
signal $\varepsilon \in \left\{ -1,1\right\} $ in the Lorentzian $3$-space,
denote by $\mathcal{Q}_{L}\left( S\right) $ the set of all polynomials $Q\in 
\mathbb{R}\left[ x,y\right] $ verifying $Q\left( K,H\right) \equiv 0$.

\begin{enumerate}
\item[i.] If $S$ is a right cylinder, then $\mathcal{Q}_{L}\left( S\right)
=\left\{ Q\in \mathbb{R}\left[ x,y\right] \,;\,r\in Rad_{L,\varepsilon
}\left( Q\right) \right\} $;

\item[ii.] If $S$ is not a right cylinder, then $\mathcal{Q}_{L}\left(
S\right) $ is the ideal in $\mathbb{R}\left[ x,y\right] $ generated by the
polynomial $xr^{2}-2ry+\varepsilon $.
\end{enumerate}
\end{theorem}

\begin{proof}
Suppose first that $S$ is a right cylinder. It is well known that the
curvatures of $S$ are%
\begin{equation*}
K\equiv 0\qquad \text{and}\qquad H=\tfrac{\varepsilon }{2r}\text{.}
\end{equation*}%
Then the equality%
\begin{equation*}
\mathcal{Q}_{L}\left( S\right) =\left\{ Q\in \mathbb{R}\left[ x,y\right]
\,;\,r\in Rad_{L,\varepsilon }\left( Q\right) \right\} 
\end{equation*}%
is an immediate consequence of the definitions of $Rad_{L,\varepsilon
}\left( Q\right) $ and $\mathcal{Q}_{L}\left( S\right) $. We observe that $%
\mathcal{Q}_{L}\left( S\right) \neq \emptyset $ because a simple calculation
shows that $xr^{2}-2ry+\varepsilon \in \mathcal{Q}_{L}\left( S\right) $.

\bigskip

Suppose now that $S$ is not a right cylinder. By Section \ref{Tubular
Surface} that the curvatures of $S$ are%
\begin{equation*}
K=\varepsilon \dfrac{\varepsilon _{B}\kappa \mu }{r\left( 1+\varepsilon
_{B}r\kappa \mu \right) }\mathbb{\qquad }\text{and}\qquad H=\varepsilon 
\dfrac{2\varepsilon _{B}r\kappa \mu +1}{2r\left( 1+\varepsilon _{B}r\kappa
\mu \right) }.
\end{equation*}%
Again a simple calculation shows that $xr^{2}-2ry+\varepsilon \in \mathcal{Q}%
_{L}\left( S\right) $ and therefore%
\begin{equation*}
\left\langle xr^{2}-2ry+\varepsilon \right\rangle \subset \mathcal{Q}%
_{L}\left( S\right) \text{.}
\end{equation*}%
Now, take a polynomial $Q$ in $\mathcal{Q}\left( S\right) $. From Corolary %
\ref{Teo_Factorization_epsilon}, to prove that $Q$ is in the ideal in $%
\mathbb{R}\left[ x,y\right] $ generated by $xr^{2}-2ry+\varepsilon $, it
sufficies to verify that%
\begin{equation*}
Q\left( \frac{\varepsilon x}{r},\frac{\varepsilon \left( xr+1\right) }{2r}%
\right) \equiv 0\text{.}
\end{equation*}%
However, if the polynomial $Q\left( \frac{\varepsilon x}{r},\frac{%
\varepsilon \left( xr+1\right) }{2r}\right) $ is not identically null, we
can repeat the arguments used in Theorem \ref{Teo S(Q) Lorentzian} to
conclude that $S$ must be a right cylinder. This contradiction show $Q$
belongs to the ideal generated by $xr^{2}-2ry+\varepsilon $ and consequently
that%
\begin{equation*}
\mathcal{Q}_{L}\left( S\right) \subset \left\langle xr^{2}-2ry+\varepsilon
\right\rangle \text{.}
\end{equation*}

\hfill \hfill \hfill
\end{proof}

\begin{corollary}[Lorentzian Linear problem]
Every non degenerated regular tubular surface in the Lorentzian $3$-space is
a linear Weingarten surface. More precisely, their curvatures verify%
\begin{equation*}
Kr^{2}-2rH+\varepsilon =0\text{,}
\end{equation*}%
where $r>0$ and $\varepsilon \in \left\{ -1,1\right\} $ are their radii and
signals.
\end{corollary}

\bigskip

Consider a (non degenerated in the Lorentzian case) regular surface $S$ in
the Euclidean, Lorentzian or Hyperbolic $3$-spaces. We say that a nonlinear
polynomial $Q\in \mathcal{Q}\left( S\right) $ defines a \textbf{true
nonlinear relation} (between the curvatures of $S$) when $Q$ does not have a
non trivial divisor in $\mathcal{Q}\left( S\right) $. Of course, nonllinear
irreducible polinomials in $\mathcal{Q}\left( S\right) $ (if they existed)
give rise to true nonlinear relations, however true nonlinear relations can
be obtained by reduclible polynomials.

\begin{corollary}[Lorentzian Nonlinear problem]
The right cylinders are the only non degenerated regular tubular surfaces in
the Lorentzian $3$-space whose curvatures vanishes a true nonlinear
polynomial relation.
\end{corollary}

\bigskip

For completeness, we present below the hyperbolic versions of the results
above:

\begin{theorem}[Hyperbolic $\mathcal{Q}\left( S\right) $ problem]
Given a regular tubular surface $S$ of radius $r>0$ in the Hyperbolic $3$%
-space, denote by $\mathcal{Q}_{H}\left( S\right) $ the set of all
polynomials $Q\in \mathbb{R}\left[ x,y\right] $ verifying $Q\left(
K,H\right) \equiv 0$.

\begin{enumerate}
\item[i.] If $S$ is a right cylinder, then $\mathcal{Q}_{H}\left( S\right)
=\left\{ Q\in \mathbb{R}\left[ x,y\right] \,;\,r\in Rad_{L}\left( Q\right)
\right\} $;

\item[ii.] If $S$ is not a right cylinder, then $\mathcal{Q}_{H}\left(
S\right) $ is the ideal in $\mathbb{R}\left[ x,y\right] $ generated by the
polynomial $x\sinh ^{2}r-2y\left( \sinh r\right) +1$.
\end{enumerate}
\end{theorem}

\begin{corollary}[Hyperbolic Linear problem]
Every regular tubular surface in the Hyperbolic $3$-space is a linear
Weingarten surface. More precisely, their curvatures verify%
\begin{equation*}
K\sinh ^{2}r-2H\left( \sinh r\right) +1=0\text{,}
\end{equation*}%
where $r>0$ are their radii.
\end{corollary}

\begin{corollary}[Hyperbolic Nonlinear problem]
The right cylinders are the only regular tubular surfaces in the Hyperbolic $%
3$-space whose curvatures vanishes a true nonlinear polynomial relation.
\end{corollary}

\bigskip

To finish this section, we observe that the algebraic approach adopted in
this article can be also used to classify Weingarten tubular surfaces
verfying a polynomial relation between its principal curvatures. We present
below only the Euclidean versions of the theorems and we leave to the reader
the task of stating its Lorentzian and Hyperbolic versions.

\begin{theorem}[Euclidean $\mathcal{S}^{\prime }\left( Q\right) $ problem]
Given a polynomial $Q\left( x,y\right) \in \mathbb{R}\left[ x,y\right] $,
denote by $\mathcal{S}_{E}^{\prime }\left( Q\right) $ the set of all regular
tubular surfaces in Euclidean $3$-space whose principal curvatures $k_{1}$ $%
\leq k_{2}$ verify $Q\left( k_{1},k_{2}\right) \equiv 0$. Then, the elements
of $\mathcal{S}_{E}^{\prime }\left( Q\right) $ are:

\begin{enumerate}
\item[i.] The right cylinders with radius in the set (called \textbf{%
Euclidean principal radius} of $Q$)%
\begin{equation*}
Rad_{E}^{\prime }\left( Q\right) =\left\{ r\in \left( 0,+\infty \right) 
\text{ };\text{ }Q\left( 0,\tfrac{1}{r}\right) =0\right\} \text{;}
\end{equation*}

\item[ii.] The regular tubular surfaces with radius in the set (called 
\textbf{Euclidean principal star radius} of $Q$)%
\begin{equation*}
\left( Rad_{E}^{\prime }\right) ^{\ast }\left( Q\right) =\left\{ r\in
Rad_{E}\left( Q\right) \,;\,Q\left( x,y\right) \in \left\langle y-\frac{1}{r}%
\right\rangle \right\} \text{.}
\end{equation*}
\end{enumerate}
\end{theorem}

\begin{corollary}[Euclidean]
Consider a polynomial $Q\left( x,y\right) \in \mathbb{R}\left[ x,y\right] $.
There is a regular tubular surfaces in Euclidean $3$-space verifying $%
Q\left( k_{1},k_{2}\right) \equiv 0$ if and only if $Rad_{E}^{\prime }\left(
Q\right) \neq \emptyset $. Moreover, the radii of all such surfaces belongs
to $\left( Rad_{E}^{\prime }\right) ^{\ast }\left( Q\right) $.
\end{corollary}

\begin{theorem}[Euclidean $\mathcal{Q}\left( S\right) $ problem]
Consider a regular tubular surface $S$ of radius $r>0$ in Euclidean $3$%
-space and let $k_{1}\leq k_{2}$ be its principal curvatures. Denote by $%
\mathcal{Q}_{E}^{\prime }\left( S\right) $ the set of all polynomials $Q\in 
\mathbb{R}\left[ x,y\right] $ verifying $Q\left( K,H\right) \equiv 0$.

\begin{enumerate}
\item[i.] If $S$ is a right cylinder, then $\mathcal{Q}_{E}^{\prime }\left(
S\right) =\left\{ Q\in \mathbb{R}\left[ x,y\right] \,;\,r\in Rad_{E}^{\prime
}\left( Q\right) \right\} $;

\item[ii.] If $S$ is not a right cylinder, then $\mathcal{Q}_{E}^{\prime
}\left( S\right) $ is the ideal in $\mathbb{R}\left[ x,y\right] $ generated
by the polynomial $y-\frac{1}{r}$.
\end{enumerate}
\end{theorem}

\bigskip

Since $k_{2}$ is constante, the linear problem for this kind of Weingarten
surfaces is trivial (in fact, the principal curvatures of any regular
tubular surface vanishes the polynomial $y-\frac{1}{r}$). However, the
nonlinear problem remains an interesting corollary of the above theorem (the
definition of true nonlinear polynomial relations in this context is
evident).

\begin{corollary}[Euclidean Nonlinear problem]
\label{Nonlinear intro}The right cylinders are the only tubular surfaces in
Euclidean $3$-space whose principal curvatures vanishes true nonlinear
polynomial relations.
\end{corollary}

\bigskip

Again we can state many particular results from the previous results. For
example:

\begin{corollary}[Euclidean]
Given real numbers $a,b,c$ such that $\left( a,b\right) \neq \left(
0,0\right) $, consider the polynomial%
\begin{equation*}
Q\left( x,y\right) =ax+by-c\in \mathbb{R}\left[ x,y\right] \text{.}
\end{equation*}%
Then, $\mathcal{S}_{E}^{\prime }\left( Q\right) \neq \emptyset $ if and only
if $b,c=0$ or $bc>0$. More precisely, we have:

\begin{enumerate}
\item[i.] If $b,c=0$, then $\mathcal{S}_{E}^{\prime }\left( Q\right) $ is
the set of all right cylinders of any radius;

\item[ii.] If $bc>0$ and $a=0$, then $\mathcal{S}_{E}^{\prime }\left(
Q\right) $ is the set of all tubular surfaces of radius $\frac{b}{c}$;

\item[iii.] If $bc>0$ and $a\neq 0$, then $\mathcal{S}_{E}^{\prime }\left(
Q\right) $ is the set of all right cylinders of radius $\frac{b}{c}$.
\end{enumerate}
\end{corollary}

\bigskip

The authors would like to express his sinceres gratitude to professor
Francisco Fontenele for the carefull reading of the initial versions of this
text and for your helpfull sugestions.

\end{document}